\documentclass[preprint]{elsarticle}

\usepackage{amssymb}

\usepackage{epsfig}

\newtheorem{thm}{Theorem}
\newtheorem{lem}{Lemma}
\newtheorem{cor}{Corollary}
\newdefinition{rmk}{Remark}
\newtheorem{prop}{Proposition}
\newcommand{\finish}{\hfill$\Box$\vspace{0.2cm}}
\newcommand{\prf}{\noindent{\bf Proof:\ }}

\newcommand{\E}{{\rm I \!E}}
\newcommand{\p}{{\rm I \!P}}

\begin{document}

\title{Stock Price Processes with Infinite Source Poisson Agents}

\author{ M. \c{C}a\u{g}lar}\fnref{fn1}
\ead{mcaglar@ku.edu.tr}
\address{Department of Mathematics, Ko\c{c} University,
 Sariyer, 34450 Istanbul, Turkey}
\fntext[fn1]{Supported by TUBITAK Project No. 109T665 Novel
Stochastic Processes for Stock Prices and their Limits}

\begin{abstract}
We construct a general stochastic process and prove weak convergence
results. It is scaled in space and through the parameters of its
distribution. We show that our simplified scaling is equivalent to
time scaling used frequently. The process is constructed as an
integral with respect to a Poisson random measure which governs
several parameters of  trading agents in the context of stock
prices.  When the trading occurs more frequently and in smaller
quantities, the limit is a fractional Brownian motion. In contrast,
a stable L\'{e}vy  motion is obtained if the rate of trading
decreases while its effect rate increases.
\end{abstract}


\maketitle


\section{Introduction}

Weak convergence  of scaled input processes has been studied
extensively over the last decade
\cite{bayrak,caglar,fasen,gaigalas,jedidi, kaj,scharf,klup,
konstan,mikres,mikosch, pipiras,resnick}. The limit is a fractional
Brownian motion (fBm) or a L\'evy process depending on the
particular scaling. While the motivation of such analysis originates
from data traffic in telecommunications, both fBm and L\'evy
processes have recently become prevalent in finance. Thereby, we
construct a general stochastic process based on a Poisson random
measure $N$, interpret it as stock price process and prove weak
convergence results.

 We consider a
 real valued process of the form
\begin{equation}  \label{1}
Z_n(t)=\int_{-\infty}^{\infty}\int_0^{\infty}\int_{-\infty}^{\infty}
 {a_n}\, r\, u\left[f\left(\frac{t-s}{u}\right)-
f\left(\frac{-s}{u}\right)\right] N_n(ds,du,dr)
\end{equation}
where $f$ is a deterministic function satisfying Lipschitz
condition, $(a_n)$ is a scaling sequence, and $r$ and $u$ are marks
of the resulting Poisson point process with  $s$ denoting the time.
The mean measure on $U$ is either a probability measure or an
infinite measure on $\mathbb{R}_+$. The process $Z_n$
 depends on the scaling parameter $n\in \mathbb{R}_+$ through not
 only $a_n$, but also through the mean
measure $\mu_n$ of $N_n$ which is taken to have a regularly
varying form in $u$ to comply with the long-range dependence
property of teletraffic or financial data. After centering the
process $Z_n$, we can obtain either an fBm or a stable L\'evy
motion depending on the particular scaling of the mean measure of
$N_n$ and the factor $a_n$ as $n\rightarrow \infty$. While fBm is
a self-similar and long-range dependent model, a L\'{e}vy process
has independent increments and self-similarity exists without
long-range dependence.

Our main contribution is the generalization of the previous results
\cite{kaj,mikres} with a specific linear form  for  $f$ and
\cite{konstan} with an increasing $f$ satisfying some other
technical conditions,  to those with Lipschitz functions. In this
case, not only the proofs require more work, but we need Lipschitz
assumptions on the derivative of $f$ as well. We
 also show that the time scaling used in previous work
 can be replaced by parameter scaling
 of the  distributions of the relevant random
 variables. The time scaling has been interpreted as `birds-eye' description
 of a process, which is not necessary when the scaling is interpreted
 in terms of its parameters.  Inspired by \cite{mandel}, we unify the results for
 general forms of $f$ with less stringent conditions in some
 cases.  In \cite{mandel},  it  is  noted that a fractional
 Brownian  motion  with $H>1/2$  can be approximated if the pulse $f$
 is continuous on
$\mathbb{R}$ and has a compact support. As an alternative extension,
we consider continuous $f$ with no compact support while
constructing (\ref{1}) with $\tilde{N}_n=N_n-\mu_n$ as a centered
process. A secondary generalization of previous work is the
consideration of $Z$ as a more general process than workload which
would be positive by definition. We allow a signed process through
the choice of real valued rate $r$.

In related work \cite{fasen,mikosch}, the Poisson random measure
is replaced by a general arrival process and a cluster Poisson
process, respectively. Ergodicity is required for the limit
theorems in the general arrival  case as well. On the other hand,
most of the previous studies on scaled input processes are named
infinite source Poisson models due to the assumption of Poisson
arrivals.

As for application in finance, the process $Z_n$ can be
interpreted as the price of a stock. The interpretation of
(\ref{1}) as a stock price process has been first presented in
\cite{akcay,world}. Our aim is to construct a model
 involving the behavior of agents
that can be parameterized and estimated from data, yet having
well-known  stochastic processes  as its limits. While the
limiting models fit well to financial data, they do not involve
the physical parameters of the trading agents. Agent based
modeling is widely used to find a model that best fits stock price
processes. In some studies, agents are divided into two groups,
mostly named as chartists and fundamentalists. In these studies,
the two agent groups have different demand functions for the stock
\cite{chartist-demand1,chartist-demand2,chartist-demand3} and the
price is generally determined via the total excess demand
\cite{demand4,demand5,demand1,demand3,demand2}. In (\ref{1}), the
arrival time $s$, the rate $r$ and the duration $u$ of the effect
of an order are all governed by the Poisson random measure. Under
the assumption of positive correlation
 between the total net demand and the price change, we expect that a
 buy order of
an agent increases the price whereas a sell order decreases it.
Each order has an effect  proportional to its volume and duration.
The duration of the effect is assumed to follow a heavy tailed
distribution. This effect starts when the order is given,
increases to a maximum which is proportional to the  total order
amount, and then starts decreasing until it vanishes after a
finite time. Alternatively, its effect may last for some time and
leave the price at a changed level on and after the time of
transaction. The logarithm of the stock price is found by
aggregating the incremental effects of orders placed by all active
agents in $[0,t]$. As a semi-martingale, our process does not
allow for arbitrage. It is a novel model which is alternative to
existing agent based constructions that are semi-Markov processes
\cite{bayrak}.

We prove weak convergence results   for the two different measures
used for the duration, separately. As for their connection, the
infinite measure being the limit of the probability measure is
proved in Theorem 1. Although a probability measure is appropriate
for the duration, its limiting form  is able to capture the
essential tail behavior which represents long-range dependence and
self-similarity properties. Moreover, the scalings with the
limiting measure is simpler to interpret as follows. An fBm limit
is found as the frequency of trading increases and the price
effect of the orders decreases as shown in Theorem 3. In contrast,
a stable process which is a scaled version of a stable L\'{e}vy
motion is obtained in Theorem 5 when the rate of trading decreases
and the effect of the orders increases. Theorems 2 and 4 are
concerned with fBm and stable limits, respectively, for the case
of probability measure for the duration. The hypotheses in these
theorems involve extra scaling needed for obtaining the limiting
measure in addition to the scaling of the rates and effects. The
stable process obtained in the limit has a skewness parameter that
depends on the distribution of the rate $r$ as it can take
positive or negative values.

 The paper is organized as follows. In Section 2, we define the
 ingredients of the workload model, namely the random variables
 associated with the process to be constructed.
  The price process is defined in Section 3 with various
 assumptions on the parameters of the Poisson random measure.
  Section 4 includes  limit theorems for fractional Brownian motion.
   Finally, limit theorems for stable L\'{e}vy
 motion are proved in Section 5.


\section{Infinite Source Poisson Agents}

 In this section, we state our main assumptions and notation for
 constructing a stock price process.  We assume potentially an
 infinite pool of agents and gather the agents' trading processes
 under a Poisson random measure. Each agent's  trading starts
 according to the underlying Poisson process and it ends after
 a random amount of time. Although the same agents may be returning
 for further transaction, new arrivals are assumed to be
 independent and identically distributed.

 The agents are called \emph{infinite source Poisson} as a term
 borrowed from traffic modeling in telecommunications based on
 Poisson random measures (e.g. \cite{caglar,mikres,mikosch}).
 Since self-similarity and long-range dependence are common
 statistical properties observed in both Internet traffic and financial data,
  the stochastic models also have similar features.

In \cite{bayrak}, a finite number of identically distributed
semi-Markov processes representing the trading states of the agents
over time are aggregated to form the price process. The number of
semi-Markov processes, equivalently, the number of agents has been
taken to infinity only in the limit. On the other hand, our price
model is a stationary process  with an infinite source of agents
that arrive according to a Poisson process. We concentrate on
scaling the other parameters that have physical interpretations for
obtaining the limiting stochastic processes.

Let $(\Omega, \mathcal{F},\mathbb{P})$ be a probability space. Let
$\mathcal{B}_{\mathbb{R}}$ denote the Borel $\sigma$-algebra on
$\mathbb{R}$. Let $N$ be a Poisson random measure on
$(\mathbb{R}\times\mathbb{R}_+\times\mathbb{R},\mathcal{B}_\mathbb{R}\otimes
\mathcal{B}_{\mathbb{R}_+}\otimes\mathcal{B}_{\mathbb{R}})$ with
mean measure
\begin{equation}
\mu(ds,du,dr)=\lambda  ds\, \nu(du)\, \gamma(dr) \label{mean}
\end{equation}
 where
$\lambda>0$ is the arrival rate of the underlying Poisson process,
$1<\delta<2$, $\gamma$ is the distribution of a random variable $R$
and
 $\nu$ is either a probability measure that satisfies
\begin{equation}
 \int_{u}^{\infty} \nu(dy) \sim h(u)  \frac{u^{-\delta}}{\delta}
\;\;\; \mbox{as} \;\;\; u\rightarrow \infty   \label{nu}
\end{equation}
where $h$ is a slowly varying function at infinity, that is, $h$ is
such that
\begin{equation}  \label{slow}
\lim_{n\rightarrow \infty}h(un)/h(n)=1
\end{equation}
 or a measure  given by
\begin{equation}
\nu(du) = u^{-\delta-1} \, du \label{mandelu} \: .
\end{equation}
 Each atom $(S_j,U_j,R_j)$ of $N$ can be interpreted as an order
from an agent, where $S_j$ is the arrival time of the order, $U_j$
is the duration of its effect on the price, and $R_j$ denotes its
rate which also plays the role of conversion to monetary units.
The sign of $R_j$ could be positive or negative depending on the
order being a buy or sell order, respectively. Under assumption
(\ref{mandelu}) for the measure $\nu$, the duration $U$ is
obtained from a diffuse measure on $(0,\infty)$ with no
probability distribution. This case is studied for suppressing the
less significant details in the proofs of the convergence
theorems. In case of (\ref{nu}), $U$ follows a heavy tailed
distribution with finite mean but infinite variance, and the
convergence proofs involve the function $h$. Although the latter
case is physically more meaningful, the scalings are  more
involved as well in the limit theorems for the self-similar
processes fBm and stable L\'{e}vy motion.

 Let
$K:\mathbb{R}\times\mathbb{R}_+\times\mathbb{R}\rightarrow\mathbb{R}$
denote an effect function such that $K(x,\cdot,\cdot)\equiv 0$ if
$x<0$. The effect of an order starting at $s$ and ending at $s+u$
depends on the rate $r$ of the effect and equals $K(t-s,u,r)$ at
time $t$ provided that $s\leq t$. The function $K$ can also be
interpreted as the local dynamics of a transaction as a result
 of a buy or sell order. The rate $r$ will be connected to the quantity
 of the order which  will be elaborated further in the sequel. We specify
 a general effect function $f: \mathbb{R} \rightarrow \mathbb{R}$ that
 determines $K$ by
\begin{equation}
K(t-s,u,r)=r u\, f\left( \frac{t-s}{u}\right)  \label{K}
\end{equation}
for $t\geq s $, and $f(x)=0$ for $x<0$. We consider $f$ to be the
deterministic shape function, or pulse, for the effect which is then
shifted to the starting position $s$, scaled and amplified for the
duration $u$, and adjusted once more with the rate/conversion factor
$r$, all randomized by $N$.

The price process will be constructed as a sum of randomized pulses.
Characterization of $f$ that will yield an fBm or a L\'{e}vy motion
for the price process is of interest, provided that the parameters
in (\ref{mean}) are appropriately scaled. In the workload processes
 studied in \cite{kaj,konstan,mikres,mikosch}, $f$ has the
following form
\begin{equation}  \label{increasingf}
f(x)=\left\{ \begin{array}{ll} x\wedge 1 & \qquad x\geq 0 \\
       0 & \qquad  x<0 \end{array} \right.
\end{equation}
which represents an increasing input, or effect, with unit rate on
$[0,1]$ and remains constant thereafter. In \cite{mandel}, the pulse
\begin{equation}  \label{mandelf}
f(x)=\left\{ \begin{array}{ll} 1/2-|x-1/2|  & \qquad 0\leq x \leq 1 \\
       0 & \qquad  \mbox{otherwise} \end{array} \right.
\end{equation}
is considered for the aim of approximating an fBm. It has a compact
support representing a limited effect that vanishes after the
duration of the pulse. These special pulses, in other words effect
functions, are sketched in  Fig.1. General Lipschitz functions are
also considered in \cite{mandel}.

\begin{figure}
\centering \epsfig{file=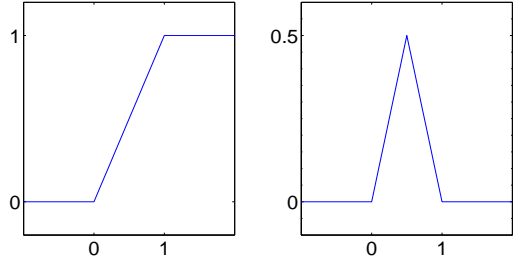, height=1.5in, width=3.3in}
\begin{center} \caption{Sample pulses of different kind}
\end{center}
\end{figure}

 The effect function is left unspecified in
\cite{jedidi,scharf,klup} but with conditions on the tail properties
of its distribution for large times. The duration is not
parameterized in contrast to the present work. In \cite{klup}, the
effect may last indefinitely although it decreases in time with a
regularly varying tail. In \cite{jedidi}, each effect is assumed to
converge for large times to a finite random variable which has a
distribution with a regularly varying tail. The stable limits are
outlined in \cite{scharf} with several interesting special cases. In
\cite{konstan}, the effect function is deterministic which is
randomized through a random variable for duration as in our case,
but with an additional assumption that the effect function itself
also has a regularly varying tail.
 On the other hand, random effect functions have
 been also considered also in \cite{kurtz} where central limit theorems are
 proved under general conditions. As a special case,
 the effect function could be a compound Poisson
process as in \cite{caglar,kaj}. This could be used to model the buy
or sell transactions in smaller quantities for a given order in the
present work. A semimartingale is assumed for the rate of the effect
in \cite{bayrak}.  In applications, $K$ can be estimated to match
the local dynamics of the price change or the workload.


\section{Price Process}

Let the price process be given by $Y=\{Y(t): Y(t)=\exp Z(t), \, t
\in \mathbb{R}_+\}$ where $Z$ is the log-price process to be
constructed in this section.    We aim to introduce a stochastic
process which is sufficiently general to approximate an fBm or a
L\'{e}vy motion, and has an adequate number of physical parameters
that can be estimated from data. The effect function and the Poisson
random measure described in Section 2 will be the main ingredients.
\begin{rmk}
Previous models that involve heterogeneous agents usually classify
them into two separate groups as chartists and fundamentalists
according to their trading behavior \cite{bayrak}. This can
clearly be generalized to several types of agents.
  The total effect $Z^i$ from agents of type $i$, $i=1,\ldots,I$,
  can be  further aggregated to form $Z$ as
\[
Z(t)=\sum_{i=1}^n Z^i(t) \;\;\;\;\;\;\;\;\;\; t\geq 0.
\]
Then, the price $Y$ at time $t$ is given by  $Y(t)=e^{Z(t)}$ as
before.
\end{rmk}

    We form the log-price process $Z$ by aggregating
    the randomized effects. More precisely,
the difference in the effect amplitudes at times $0$ and $t$ are
integrated with respect to the Poisson random measure $N$ to yield
\[
Z(t)=\int_{-\infty}^{t}\int_0^{\infty}\int_{-\infty}^{\infty}
 (K (t-s,u,r)-K (-s,u,r))N(ds,du,dr)\qquad t\geq 0.
 \]
Since the underlying Poisson process has been going on long before
time 0, $Z$ has stationary increments and  $Z (0)=0$ by
construction. We think $Z(t)-Z (0)$ as the sum of all effects due
to all active agents between times 0 and $t$. We assume that $K$
has the form (\ref{K}) as before, and write
\begin{equation}  \label{Zf}
Z(t) = \int_{-\infty}^{\infty} \int_0^\infty
\int_{-\infty}^{\infty}ru\,\left[f\left(\frac{t-s}{u}
  \right)
  -f\left(\frac{-s}{u}\right)\right]\, N(ds,du,dr).
  \end{equation}
The  propositions below give the sufficient conditions for $Z$ to be
well-defined for a finite measure and a specific $\sigma$-finite
measure $\nu$ in Equation (\ref{mean}), respectively. Let $\phi$
denote the characteristic function of $Z$, that is, $\phi(\xi)=
\mathbb{E} e^{i\xi Z}$ for $\xi \in \mathbb{R}$.

\begin{prop}  \label{prop1}
 Suppose that  $\nu$ is a probability measure
 satisfying $\int u\nu(du)<\infty$, $\mathbb{E}|R|<\infty$
 and $f$ is a Lipschitz continuous function on $\mathbb{R}$ with
 $f(x)=0$ for all $x\leq 0$   and $f(x)=f(1)$ for all $x\geq 1$. Then,
$Z(t)$, $t\geq 0$, is a finite random variable a.s. with
characteristic function
  \[
 \phi(\xi)=\exp \int_{\mathbb{R}} \int_{\mathbb{R}_+}
  \int_{\mathbb{R}} \left\{\exp \left[i\xi \,ru\left(f\left(\frac{t-s}{u}
  \right)
  -f\left(\frac{-s}{u}\right)\right)\right] -1 \right\}
    \lambda \, ds\, \nu(du)\, \gamma(dr)
  \]
\end{prop}

\noindent{\bf Proof:} The integral $\int h d\mu$ of a deterministic
function $h$ with respect to a Poisson random measure defines a
finite random variable if $\int |h \wedge 1| \, d\mu < \infty$ where
$\mu$ is the mean measure. This is clearly satisfied if $\int |h| \,
d\mu < \infty$ which also implies that the random variable has a
finite expectation \cite{kingman}. Therefore, it is sufficient to
show that the expression
\[
  I  :=
\int_{-\infty}^{\infty} \int_0^{\infty}u\left|
f\left(\frac{t-s}{u}\right)- f\left(\frac{-s}{u}\right)\right|
\nu(du)ds
\]
is finite for $Z$ to be well defined as $\mathbb{E}|R|<\infty$. Note
that
\[
 I =  \int_{-\infty}^0\int_0^{\infty}u
\left|f\left(\frac{t-s}{u}
  \right)
  -f\left(\frac{-s}{u}\right)\right|\, \nu(du) \,ds +
  \int_{0}^t\int_0^{\infty}u\left|f\left(\frac{t-s}{u}
  \right)\right|\, \nu(du) \,ds
\]
 Considering the two different
regions $\{u: 0<-s/u<1<(t-s)/u\}$ and $\{u: 0<-s/u<(t-s)/u<1\}$ for
the first integral, and the regions
$\{u: (t-s)/u>1\}$ and $\{u: (t-s)/u<1\}$ for the second integral in
$I$, we get
\begin{eqnarray*}
I & = & \int_{-\infty}^0\left[ \int_{-s}^{t-s}u\left|f(1)
  -f\left(\frac{-s}{u}\right)\right|  \nu(du)   +
  \int_{t-s}^{\infty}u\,\left|f\left(\frac{t-s}{u}
  \right)
  -f\left(\frac{-s}{u}\right)\right|  \nu(du) \right] ds \\
  &  &  + \int_{0}^t \left[ \int_0^{t-s}u\,\left|f(1)\right|
 \, \nu(du) \,ds +
 \int_{t-s}^{\infty}u\,\left|f\left(\frac{t-s}{u}\right)
 \right| \,\nu(du) \right] \,ds
\end{eqnarray*}
Due to the Lipschitz hypothesis on $f$ and since $f(0)=0$,  we  have
\begin{eqnarray}
I&\leq &\int_{-\infty}^0\left[ \int_{-s}^{t-s} u \,  M
  \left(1+\frac{s}{u}\right)  \nu(du)   +
  \int_{t-s}^{\infty}u \, M  \, \frac{t }{u}\,
  \nu(du) \right] ds   \label{above} \\
 & &   + \int_{0}^t \left[ \int_0^{t-s}u\,M
 \, \nu(du)   +
 \int_{t-s}^{\infty}u\,M \frac{t-s}{u} \,\nu(du) \right] \,ds
 \nonumber
\end{eqnarray}
where $M>0$. We apply integration by parts for the inner integrals
above. For $-\infty<s<0$, integration by parts yields
\[
\int_{-s}^{t-s}(s+u)\nu(du)=tF(t-s)-\int_{-s}^{t-s}F(u)du
\]
 and we have $t\int_{t-s}^\infty
\nu(du)=t\bar{F}(t-s)$ where $F$ denotes the cumulative distribution
function (cdf) of  $U$ and $\bar{F}=1-F$. For $0<s<t$, we get
\[
\int_{0}^{t-s}u\, \nu(du)=(t-s)F(t-s)-\int_0^{t-s}F(u)du
\]
by integration by parts and we have $(t-s)\int_{t-s}^{\infty}
\nu(du)=(t-s)\bar{F}(t-s)$. Putting all expressions together and
using $\bar{F}=1-F$, we simplify (\ref{above}) as
\[
I\leq M \int_{-\infty}^0\int_{-s}^{t-s}\bar{F}(u)\, du ds+M\int_0^t
\int_0^{t-s} \bar{F}(u) \,du ds
\]
Changing the order of integration, we get
\[
  I \leq M \int_0^\infty \int_{-u}^{t-u} \bar{F}(u) \, ds du
  =Mt\mathbb{E}U
\]
which is finite by hypothesis as $\mathbb{E}U = \int u \nu(du)$. The
characteristic function of $Z$ can be found immediately from
formulae for integrals with respect to a Poisson random measure
\cite{kingman}. \finish

The  function $f$ characterized in Proposition \ref{prop1}
represents the local dynamics due to the effect of an individual
buy or sell order. The change in price, which can be
non-monotonic, occurs over a finite time $U$ and remains at the
same level thereafter. The specific shape of $f$ is left
unspecified, and so is its sign. In general, we expect a buy order
to increase the price and a sell order to decrease it. Therefore,
if $f$ is chosen to be an increasing function, then we could have
$R>0$ for a buy order and $R<0$ for a sell order. However, a
general form is assumed to leave room for modeling purposes in
view of real data and to provide mathematical generality.  The
special case (\ref{increasingf}) used in \cite{kaj,mikosch,mikres}
is a linearly increasing pulse as stated in the following
corollary.

\begin{cor} Suppose that  $\nu$ is a probability measure
 satisfying $\int u\, \nu(du)<\infty$, $\mathbb{E}|R|<\infty$ and
 $f(x)=x\wedge 1$, $x\geq 0$. Then, $Z(t)$, $t\geq 0$, is finite a.s. and
 its characteristic function is given by
\[
 \exp \int_{\mathbb{R}} \int_{\mathbb{R}_+}
  \int_{\mathbb{R}} \left[e^{  i\xi \,r \left( u \wedge (t-s)^+
  - u \wedge (-s)^+ \right)  } -1 \right]
    \lambda \, ds\, \nu(du)\, \gamma(dr)
 \]
for $\xi \in {\mathbb{R}}$.
\end{cor}

\begin{rmk}
  The log-price process $Z$ which is a semimartingale in
general,
   becomes a martingale if its mean is zero. This would be
 satisfied if $\mathbb{E} R =0$, which corresponds to symmetric
 effects from buy and sell orders, for example.
 \end{rmk}

 The following proposition is based on the results of \cite{mandel}.
 The support of $f$ is chosen as [0,1] for simplicity without loss of generalization.

\begin{prop} \label{prop2} Suppose that $\nu(du) = u^{-\delta-1} \, du$,
$\mathbb{E}R^2<\infty$ and $f:\mathbb{R}\rightarrow \mathbb{R}$ is
Lipschitz continuous with compact support $[0,1]$. Then,  $Z$ is a
well-defined random variable  with characteristic function
\[
 \exp  \int_{\mathbb{R}}\int_{\mathbb{R}_{+}}
 \int_{\mathbb{R}}\!\! \left\{e^{i\xi r u
\left[f\left(\frac{t-s}{u}\right)-f\left(\frac{-s}{u}\right)\right]}-1
 -i\xi  r u
\left[f\scriptstyle{\left(\!\frac{t-s}{u} \! \right)} -
\displaystyle{f} \scriptstyle{\left( \! \frac{-s}{u}\! \right)}
\right] \! \right\} u^{-\delta -1} \! \lambda \,ds \, du \gamma(dr)
\]
\end{prop}
\noindent \textbf{Proof:}  The Lipschitz assumption  on $f$ and that
it has compact support $[0,1]$ imply $f(0)=f(1)=0$ and
\begin{equation}  \label{varfin}
\int_{\mathbb{R}}\int_{\mathbb{R}_{+}}
\left[f\scriptstyle{\left(\!\frac{t-s}{u} \! \right)} -
\displaystyle{f} \scriptstyle{\left( \! \frac{-s}{u}\! \right)}
\right]^2u^{1-\delta}\, du \, ds < \infty
\end{equation}
for each $t>0$, and $1<\delta<2$ in particular, by
\cite[Prop.3.1]{mandel}. We sketch the usual proof for defining $Z$
as an almost sure limit of zero mean random variables, as in
\cite{mandel}. Let $ A_k=(2^{-k},2^{-k+1}],\,
k=1,2,\ldots,\,A_0=(1,\infty)$ be a partition of $\mathbb{R}_+$.
Clearly,
$$\int_{-\infty}^{\infty}\int_{A_k}\int_{-\infty}^{\infty}r u
 \left[f\scriptstyle{\left(\!\frac{t-s}{u} \! \right)} -
\displaystyle{f} \scriptstyle{\left( \! \frac{-s}{u}\! \right)}
\right]u^{-\delta -1} \lambda ds \, du \gamma(dr) =0,$$ for
$k=1,2,\ldots$, due to the form of the effect function and
boundedness of $A_k$, and also for $k=0$ in view of (\ref{varfin}).
Hence, the random variables
\begin{equation}
\int_{-\infty}^{\infty}\int_{A_k} \int_{-\infty}^{\infty}r\, u
\,\left[f\scriptstyle{\left(\!\frac{t-s}{u} \! \right)} -
\displaystyle{f} \scriptstyle{\left( \! \frac{-s}{u}\! \right)}
\right] N(ds,du,dr)\,\,\,k=0,1,\ldots, \label{mandel_random_var}
\end{equation}
are well-defined, independent and have zero expectations. One can
show that the sum of their variances given by
\begin{equation}  \label{var}
\int_{-\infty}^\infty\int_0^\infty\int_{-\infty}^\infty
   r^2 \,u^2
\left[f\scriptstyle{\left(\!\frac{t-s}{u} \! \right)} -
\displaystyle{f} \scriptstyle{\left( \! \frac{-s}{u}\! \right)}
\right] ^2 u^{-\delta -1} \lambda ds \, du \,\gamma(dr)
\end{equation}
 is finite using (\ref{varfin}) and the assumption that
 $\mathbb{E}R^2<\infty$. Therefore, the series
\[
\sum_{k=0}^\infty\int_{-\infty}^{\infty}\int_{A_k}
\int_{-\infty}^{\infty}r\, u
\left[f\scriptstyle{\left(\!\frac{t-s}{u} \! \right)} -
\displaystyle{f} \scriptstyle{\left( \! \frac{-s}{u}\! \right)}
\right]  N (ds,du,dr)\label{series}
\]
is a.s. convergent by \cite[Lemma 3.16]{kall}. The limit is denoted
by (\ref{Zf}) which is a stochastic integral in general since it may
be defined by an almost sure limit as above. Its characteristic
function at $\xi\in\mathbb{R}$ is the limit of the characteristic
functions of the partial sums of the random variables
(\ref{mandel_random_var}), since almost sure convergence implies
convergence in distribution. Now, the characteristic function of
(\ref{mandel_random_var}) can be written as
\[
 \exp  \int_{\mathbb{R}}  \int_{A_k}\! \!
\int_{\mathbb{R} }\!\! \left\{e^{i\xi r u
\left[f\left(\frac{t-s}{u}\right)-f\left(\frac{-s}{u}\right)\right]}-1
 -i\xi  r u
\left[f\scriptstyle{\left(\!\frac{t-s}{u} \! \right)} -
\displaystyle{f} \scriptstyle{\left( \! \frac{-s}{u}\! \right)}
\right]  \! \right\} u^{-\delta -1} \! \lambda \,  ds \, du
\gamma(dr)
\]
since (\ref{mandel_random_var}) has zero mean. Using the
independence of (\ref{mandel_random_var}) by disjointness of $A_k$,
$k=0,1,2,\ldots$, we can write the characteristic function of their
partial sums up to say $m\in \mathbb{Z}_+$ as
\[
 \exp  \int_{\mathbb{R}}  \int_{2^{-m}} ^\infty \!
\int_{\mathbb{R} }\!\! \left\{e^{i\xi r u
\left[f\left(\frac{t-s}{u}\right)-f\left(\frac{-s}{u}\right)\right]}-1
 -i\xi  r u \left[f\scriptstyle{\left(\!\frac{t-s}{u} \! \right)} -
\displaystyle{f} \scriptstyle{\left( \! \frac{-s}{u}\! \right)}
\right]  \! \right\} u^{-\delta -1} \! \lambda \,  ds \, du
\gamma(dr)
\]
 Due to the inequality
$|e^{ix}-1-ix|<\frac{1}{2}x^2$ for $x\in\mathbb{R}$ and by
finiteness of (\ref{var}),  dominated convergence theorem applies
and we get the result. \finish

The pulse (\ref{mandelf}) considered in \cite{mandel} is a special
case for Proposition \ref{prop2} where the process $Z$ is a
zero-mean martingale. In fact, we have $\mathbb{E}|Z|<\infty$ with
the particular pulse (\ref{mandelf}) which has the shape of an
isosceles triangle. In this case, the effect starts at time $s$ as
the buy or sell order is first given. It increases linearly as the
amount traded increases, reaches its maximum value at time
$s+\frac{u}{2}$ as the highest effect  is reached and starts
decreasing from that point on until it vanishes at time $s+u$ and
brings the price level back to the original, locally. Such functions
$f$ represent the effect of an order which takes place over a
period, and then vanishes after a while. This scenario is a milder
and time limited version of the Poisson shot-noise of \cite{klup}
where a shock in the market changes the price through a jump, but
then it tends back to its initial level by an exponential decay of
the first effect. There is no need to get back the effect as  in the
pulses of \cite{mandel} for the limit theorems. The effect function
may leave the price in a level different from the one it found at an
arrival, in which case we consider the centered process.

 We
can show the connection of the two measures in propositions
\ref{prop1} and \ref{prop2} in the next theorem.  We introduce a
scaling factor $n\in\mathbb{Z}_+$ which will be taken to infinity in
the limit. An integral with respect to a probability measure
$\nu(du)$ with regularly varying tail converges to an integral with
respect to the measure $u^{-\delta -1}du$ in the limit. We first
need the following lemma.

\begin{lem} \label{lemma1} Suppose  $f$ is a Lipschitz continuous function
 on $\mathbb{R}$ with
 $f(x)=0$ for all $x\leq 0$   and $f(x)=f(1)$ for all $x\geq 1$. Then,
\begin{equation}
\int_{-\infty}^{\infty} \int_0^{\infty} u^{1+\kappa} \left[
f\left(\frac{t-s}{u}\right)-
f\left(\frac{-s}{u}\right)\right]^{1+\kappa}
 u^{-\delta -1} \, du \, ds \: < \infty   \label{secmoment}
\end{equation}
for $1<\delta<3$ and $\kappa>0$ such that $1+\kappa >\delta$.
\end{lem}

\noindent \textbf{Proof:}  Let us call the integral
(\ref{secmoment}) as $I$. In view of the assumptions on $f$, we have
\begin{eqnarray*}
|I|\leq M^{1+\kappa} \left\{ \int_{-\infty}^0\int_{-s}^{t-s}
(u+s)^{1+\kappa} u^{-\delta -1} du \, ds+\int_{0}^t \int_{0}^{t-s}
u^{1+\kappa} u^{-\delta -1} du \,  ds \right.
\\
\;\;\;\;\;\;\;\; \left.  + \int_{0}^t
\int_{t-s}^{\infty}(t-s)^{1+\kappa} u^{-\delta-1} du \, ds
+\int_{-\infty}^0   \int_{t-s}^{\infty} t^{1+\kappa} u^{-\delta-1}
du \, ds\right\}
\end{eqnarray*}
which can be shown along the same lines as in the proof of
Proposition \ref{prop1}. Evaluating the above integrals, we find
that
\begin{eqnarray*}
\lefteqn{|I|\leq M^{1+\kappa} \,t^{2+\kappa-\delta}\left[
\frac{1}{(2 + \kappa)(2+\kappa- \delta) } + \frac{1}{(2 + \kappa)
\delta}
+\frac{1}{(2 + \kappa - \delta)(1 + \kappa - \delta) } \right.} \\
 & \qquad \qquad \qquad  \qquad \qquad \qquad  \;\;
     \qquad \qquad \qquad  \;\;\; \displaystyle{ \left.
     +\frac{1}{\delta (2 + \kappa - \delta) }
+\frac{1}{\delta(\delta -1)} \right]}
\end{eqnarray*}
\finish

\begin{thm}  \label{ICR}
Suppose that  $\nu$ is a probability measure
 satisfying $\int u \, \nu(du)<\infty$ and has a regularly varying
 tail as given in
 (\ref{nu}), the function $f:\mathbb{R}\rightarrow \mathbb{R}$ is Lipschitz
continuous with $f(x)=0$ for all $x\leq 0$, $f(x)=f(1)$ for all
$x\geq 1$ and is also differentiable with $f'$ satisfying  a
Lipschitz condition a.e., and $\mathbb{E}|R|^{1+\kappa}<\infty$
for some $0<\kappa\leq 1$ with $1+\kappa>\delta>1$. Let
\[
\mu_n (ds,du,dr)= \frac{n^{\delta}}{ h(n) }\: \lambda \, ds \,
\nu_n(du)\, \gamma(dr)
\]
where $\nu_n(du) = \nu(d(nu))$ and
\[
Z_n(t)=\int_{-\infty}^{\infty}\int_0^{\infty}
\int_{-\infty}^{\infty}   r \,u\left[f\left(\frac{t-s}{u}\right)-
f\left(\frac{-s}{u}\right)\right]N_n(ds,du,dr).
\]
Then,   $\{Z_n(t)-\mathbb{E} Z_n(t) ,t\geq 0\}$ converges in law to
 the process
\[
   \left \{ \int_{-\infty}^{\infty}\int_0^{\infty} \int_{-\infty}^{\infty}
 r \,u\left[f\left(\frac{t-s}{u}\right)-
f\left(\frac{-s}{u}\right)\right]\tilde{N}'(ds,du,dr)\; , \,t\geq
0\right\}
\]
 as \ $n\rightarrow\infty$, where $\tilde{N}'=N'-\mu'$
  for a Poisson random measure $N'$ with mean measure
  $\mu'( ds,du,dr)=  \lambda \,  u^{-\delta-1} ds \, du\,\gamma(dr)$.
\end{thm}
\prf For the convergence of finite dimensional distributions of
$\{Z_n(t)-\mathbb{E} Z_n(t) ,t\geq 0\}$, consider the
characteristic function $\mathbb{E}\exp i\sum_{k=1}^m\xi_k
[Z_n(t_k)-\mathbb{E} Z_n(t_k)]  $ for $\xi_k\in \mathbb{R}$,
$t_k\geq 0$ and $m\in \mathbb{N}$. It is given by
\begin{eqnarray}
&& \exp \int_{-\infty}^{\infty}\int_0^{\infty}
\int_{-\infty}^{\infty}\Bigg\{ e^{i\sum_{k=1}^m\xi_k r \,
u\left[f\left(\frac{t_k-s}{u}\right)-
f\left(\frac{-s}{u}\right)\right]}-1\Bigg. \nonumber\\
&&\qquad \left.-i\sum_{k=1}^m\xi_k\, r
\,u\left[f\left(\frac{t_k-s}{u}\right)-
f\left(\frac{-s}{u}\right)\right]\right\} \frac{n^{\delta}}{h(n)}
\: \lambda\,
 ds\, \nu_n( du)\,\gamma(dr) \label{chr1}.
\end{eqnarray}
We first show that the exponent in (\ref{chr1}) is bounded and then
use bounded convergence theorem to take the limit. This theorem is a
generalization of \cite[Thm.1]{kaj} with the general effect function
$f$. Although we follow the same approach as in \cite[Thm.1]{kaj},
there are more terms to bound in our case. Let
\begin{equation}  \label{g}
g(s,u,r)= e^{i\sum_{k=1}^m\xi_k r \,
u\left[f\left(\frac{t_k-s}{u}\right)-
f\left(\frac{-s}{u}\right)\right]}-1 -i\sum_{k=1}^m\xi_k \, r
\,u\left[f\left(\frac{t_k-s}{u}\right)-
f\left(\frac{-s}{u}\right)\right]\: .
\end{equation}
Using the random variable $U$, we denote the left hand side of
(\ref{nu}) as $\p\{U\geq u\}$ below. By integration by parts, the
exponent in (\ref{chr1}) is equal to
\begin{equation}
\int \int \int  \partial_u g(s,u,r) \, \p\{U>nu\}\,
\frac{n^{\delta}}{h(n)}\: \lambda \, ds \, du\, \gamma(dr)
\label{byparts}
\end{equation}
where $\partial_u$ is $\partial/\partial u$ and the hypothesis that
$\nu_n(du)=\nu(d(nu))$ is used.

\textbf{a) }Bound for the integrand of (\ref{byparts}) for large
values of $u$:

 In view of Potter bounds
\cite{bingham}, for $\epsilon>0$  there exists $n_0 \in \mathbb{N}$
such that
\[
\frac{\p\{U> nu\}}{\p\{U> n\}} \leq 2 u^{-\delta} \max
(u^{-\epsilon},u^{\epsilon})
\]
for all $n\geq n_0$ and $nu\geq n_0$, that is, $u\geq n_0/n$. Since
$\lim_{n\rightarrow \infty} \p\{U>n\}n^\delta/h(n) =C$ for some
$C>0$,  we have $\p\{U>n\}n^\delta/h(n) \leq (C+\epsilon)$ for all
$n\geq n_0'$ for some $n_0'\in \mathbb{N}$. Note that $C=1/\delta$
by (\ref{nu}). Assume $n_0'\leq n_0$ for simplicity of notation.
Therefore, we get
\begin{equation} \label{14}
\p\{U>nu \} \frac{n^{\delta}}{h(n)} \leq 2 u^{-\delta} \max
(u^{-\epsilon},u^{\epsilon})(C+\epsilon)
\end{equation}
for all $n\geq n_0$ and $u\geq n_0/n$.

In (\ref{byparts}),  we have
\begin{equation} \label{gu}
\partial_u g(s,u,r)=i\left[e^{i\sum_{k=1}^m\xi_k r \,
u\left[f\left(\frac{t_k-s}{u}\right)-
f\left(\frac{-s}{u}\right)\right]}-1 \right]\partial_u S(s,u,r)
\end{equation}
where
\begin{eqnarray*}
\partial_u S(s,u,r)& :=& \frac{\partial}{\partial u} \sum_{k=1}^m\xi_k r \,
u\left[f\left(\frac{t_k-s}{u}\right)-
f\left(\frac{-s}{u}\right)\right] \\
& = & \sum_k \xi_k r \, \left[f\left(\frac{t_k-s}{u}\right)-
f\left(\frac{-s}{u}\right)\right] \\
& & \qquad + \sum_k \xi_k r \,
\left[-f'\left(\frac{t_k-s}{u}\right)\frac{t_k-s}{u}+
f'\left(\frac{-s}{u}\right)\frac{-s}{u}\right]
\end{eqnarray*}
Now, we can bound $|\partial_u S|$ using the Lipschitz property of
$f$ and $f'$ on different regions for $u$ and $s$. Let $M>0$ and
$M'>0$ stand for the Lipschitz constants of $f$ and $f'$,
respectively, or their upper bound, whichever is larger. Let us
assume $M>M'$ for simplicity of notation.

 i) $s<0$ and $0<s+u<t_k$

Since $(t_k-s)/u>1$ and $-s/u<0$, we have
\[
\left|f \left(\frac{t_k-s}{u}\right)-f\left(\frac{-s}{u}\right)
\right|  =  \left|f(1)-f\left(\frac{-s}{u}\right) \right|  \leq
 M \left|1+\frac{s}{u}\right|
\]
 and
\begin{eqnarray*}
\left|-f'\left(\frac{t_k-s}{u}\right)\frac{t_k-s}{u}+
f'\left(\frac{-s}{u}\right)\frac{-s}{u}\right] =
\left|0-f'\left(\frac{-s}{u}\right)\frac{-s}{u}\right| \leq M'
\left| \frac{s}{u}\right|\leq M \left| \frac{s}{u}\right|
\end{eqnarray*}
due to the form of $f$ and Lipschitz assumptions. Therefore, we get
\[
|\, \partial_u S(s,u,r) | \leq  \left( M \left|1+\frac{s}{u}\right|
+ M \left| \frac{s}{u}\right| \right) \sum_k \xi_k |r|  = M\sum_k
|\xi_k|\,|r| \,
\]
since $|1+s/u|<1$ and $|s/u|<1$ in this region.

 ii) $s>0$ and $s+u<t_k$

In this region, $f'$ vanishes and $f(-s/u)=0$. Therefore, we have
\[
|\partial_u S(s,u,r) | = |f(1)| \leq M\sum_k |\xi_k|\, |r|\: .
\]

iii) $0<s<t_k$ and $t_k< s+u$

In this region, $f(-s/u)=f'(-s/u)=0$ and we get
\[
|\partial_u S(s,u,r) | \leq 2M \sum_k |\xi_k| \left|\frac{t_k-s}{u}
\right|\, |r| \: .
\]

iv) $s<0$ and $t_k<s+u$

We have
\[
\left|f \left(\frac{t_k-s}{u}\right)-f\left(\frac{-s}{u}\right)
\right| \leq M \left| \frac{t_k}{u} \right|
\]
and
\[
\left|f'
\left(\frac{t_k-s}{u}\right)\frac{t_k-s}{u}-f'\left(\frac{-s}{u}\right)\frac{-s}{u}
\right| \leq M' \left| \frac{t_k}{u} \right| + M'\left| \frac{s
t_k}{u^2}\right|\: .
\]
The corresponding bound on $|\partial_u S(s,u,r) |$ follows.

Now, we can bound the remaining terms in (\ref{gu}) by
\begin{equation}  \label{and}
  2^{1-\kappa} \sum_j |\xi_j|^\kappa |r|^\kappa u^\kappa
\left|f\left(\frac{t_j-s}{u}\right)-f\left(\frac{-s}{u}\right)
\right|^\kappa \;\;  \wedge \; \; 2
\end{equation}
using the inequalities $|e^{ix}-1| \leq 2^{1-\kappa} |x|^{\kappa}$
and $\left(\sum_j|x_j|\right)^\kappa \leq \sum_j |x_j|^\kappa$,
$0<\kappa \leq 1$ \cite{kaj}, and the fact that $|e^{ix}-1|\leq
2$. The index $k$ is replaced by $j$  in order to distinguish the
cross products of sums below. We further note that
\begin{equation}  \label{and2}
\left|f\left(\frac{t_j-s}{u}\right)-f\left(\frac{-s}{u}\right)
\right|^\kappa \leq  (2M)^\kappa  \leq 2^\kappa M
\end{equation}
since $f$ is bounded and $M^\kappa\leq M$, assuming $M>1$ for
simplicity of notation. Putting all terms together by (\ref{14}),
(\ref{and}), (\ref{and2}) and i)-iv), we find that (\ref{byparts})
is bounded as
\begin{eqnarray} \label{bound}
\lefteqn{\int  \int_{n_0/n}^\infty \int |\partial_u g(s,u,r)| \,
\p\{U>nu\}\,
\frac{n^{\delta}}{h(n)} \, \lambda \, ds \, du\, \gamma(dr)} \\
&  \leq  4 M (C+\epsilon) \displaystyle{\sum_k |\xi_k|
 \!\! \int \!\! \int_{0}^\infty  \!\!\!\! \int
|r|\, B(s,u,t_k) \max (u^{-\epsilon},u^{\epsilon})\,
 u^{-\delta}\, \lambda  \, ds \, du \, \gamma(dr) } \nonumber
\end{eqnarray}
where
\begin{eqnarray}
\lefteqn{B(s,u,t_k)= \left( 1\wedge \sum_j    M |\xi_j|^\kappa
|r|^\kappa
 u^{\kappa}1_{\{s\leq t_j\}} \right)} \nonumber \\
 & \qquad  \displaystyle{ . \: \left[  1_{R_{1,k}}+   1_{R_{2,k}}+ 2 \frac{t_k-s}{u}
 \, 1_{R_{3,k}}+  \left(2   \frac{t_k}{u}   +
\frac{|s| t_k}{u^2} \right) 1_{R_{4,k}} \right]  }\label{B}
\end{eqnarray}
and $R_{1,k},\ldots ,R_{4,k}$ denote the regions in i)-iv).
Since $[  1\wedge \sum_j    M |\xi_j|^\kappa |r|^\kappa
 u^{\kappa}1_{\{s\leq t_j\}}] $ $  \leq 1 $,  we can write
\begin{eqnarray} \label{nextB}
 B(s,u,t_k) \leq   1_{R_{1,k}}+  2 \frac{t_k-s}{u}
 \, 1_{R_{3,k}}+  \left(2   \frac{t_k}{u}   +
\frac{|s| t_k}{u^2} \right) 1_{R_{4,k}} \\ + \left( 1\wedge \sum_j M
|\xi_j|^\kappa |r|^\kappa \nonumber
 u^{\kappa}1_{\{s\leq t_j\}} \right)   1_{R_{2,k}}
\end{eqnarray}
We keep the extra bounding term for $R_{2,k}$, as the integration in
this region is more delicate. For fixed $k\in \{1,\ldots n\}$,
$R_{1,k}, \ldots,R_{4,k}$ are depicted in Fig.2. If we choose
$\epsilon>0$ such that
\begin{equation} \label{epkappa}
1<\delta -\epsilon<\delta <\delta+\epsilon<1+\kappa \: ,
\end{equation}
then the right hand side of (\ref{bound})  is  finite as as shown in
Appendix A.

\begin{figure}
\centering \epsfig{file=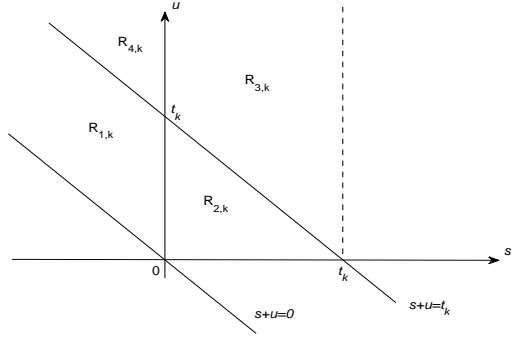, height=2in, width=3.3in}
\begin{center} \caption{Subregions considered for $(s,u)$}
\end{center}
\end{figure}

\textbf{b)} Bound for the integrand of (\ref{byparts}) for small
values of $u$:

We now consider $u\leq n_0/n \leq 1$ as $n\geq n_0$. We use Markov
inequality for $\p\{U\geq nu\}$ as in \cite{kaj}, together with the
bounds (\ref{and}), (\ref{and2}) and i)-iv) in our case. We have $
\p\{U\geq nu\} \leq  \E U / (nu)$ by Markov's inequality. Therefore,
we get
\begin{equation}  \label{bound4}
 |\partial_u g(s,u,r)| \, \p\{U>nu\}\,
\frac{n^{\delta}}{h(n)}  \leq  2M \frac{\E U}{u}
\frac{n^{\delta-1}}{h(n)} \, |r| \sum_k|\xi_k|  B(s,u,t_k)
\end{equation}
From (\ref{nextB}), we can write
\begin{eqnarray*}
B(s,u,t_k) \leq
  1_{R_{1,k}}\; +\;  2 \frac{t_k-s}{u}
 \, 1_{R_{3,k}} & \!\! \!\!\!\!+ &\!\!\!\!\!\!  \left(2   \frac{t_k}{u}   +
\frac{|s| t_k}{u^2} \right)  1_{R_4,k}    \\ & +& \!\! 1_{R_{2,k}}
\sum_j M |\xi_j|^\kappa |r|^\kappa u^{\kappa}1_{ \{s\leq t_j\}}
\end{eqnarray*}
considering that $u$ is bounded as $u\leq 1$. Now, we have
\begin{equation} \label{bound3}
1=   u^{\delta+\epsilon-1}u^{1-\delta-\epsilon} \leq
n_0^{\delta+\epsilon-1} n^{-\epsilon}n^{1-\delta} u^{1
-\delta-\epsilon} \leq n_0^{\delta+\epsilon-1} h(n) n^{1-\delta}u^{1
-\delta-\epsilon}
\end{equation}
since $u\leq n_0/n$ and $n^{-\epsilon}\leq h(n)$ for the slowly
varying function $h$ when $n$ is sufficiently large \cite{kaj}.
Using (\ref{bound3}) to increase the righthand side of
(\ref{bound4}) and in view of (\ref{nextB}), we get
\begin{eqnarray*}
 \lefteqn{ |\partial_u g(s,u,r)| \, \p\{U>nu\}\,
\frac{n^{\delta}}{h(n)}   \leq   2  M\E U
  \, n_0^{\delta+\epsilon-1}\, |r| u^{-\delta-\epsilon} \sum_k|\xi_k|    }
\\   &   .\,  \displaystyle{
\left[  1_{R_{1,k}}+  2 \frac{t_k-s}{u}
 \, 1_{R_{3,k}}+  \left(2   \frac{t_k}{u}   +
\frac{|s| t_k}{u^2} \right) 1_{R_{4,k}} + 1_{R_{2,k}} \sum_j    M
|\xi_j|^\kappa   |r|^{\kappa} u^{\kappa} 1_{\{s\leq t_j\}} \right]}
\end{eqnarray*}
which is integrable over $0<u<1$, as in part a), in view of the
computations in Appendix A.

As a result of  a) and b),  the integrand in (\ref{byparts}) is
bounded by an integrable function. Therefore, we can use
dominated convergence theorem to find that
\begin{equation} \label{kucuklimit}
\lim_{n\rightarrow \infty} \p\{U>nu\}\, \frac{n^{\delta}}{h(n)}
=\frac{u^{-\delta}}{\delta}
\end{equation}
by (\ref{nu}) and (\ref{slow}), and then revert (\ref{byparts}) by
another integration by parts to get the limit of (\ref{chr1}) as
\[
  \exp \int_{-\infty}^{\infty}\int_0^{\infty}
\int_{-\infty}^{\infty}g(s,u,r)\: \lambda u^{-\delta -1}
 ds\,   du\,\gamma(dr)
\]
where $g$ is as in (\ref{g}). It can be shown as in Proposition
\ref{prop2} that the above
 characteristic function
and the corresponding process are well defined since $|e^{ix}-1-ix|$
is bounded by $|x|$. Hence, we have shown the convergence of finite
dimensional distributions.

   To prove weak convergence in the Skorohod topology on
   $D(0,\infty)$, we first observe that
\begin{eqnarray}   \label{upper}
\lefteqn{\mathbb{E}|Z_n(t)-\mathbb{E}Z_n(t)|^{1+\kappa} \leq }\\
& \;\;\;\;\;\;\;\;\;\;\;\;\;\;\;\;\; \displaystyle{2
\mathbb{E}|R|^{1+\kappa} \int_0^\infty \int_{-\infty}^t
  u^{1+\kappa}\!\! \left[f\left(\frac{t-s}{u}\right)-
f\left(\frac{-s}{u}\right)\right]^{1+\kappa} \!\! \nonumber
\frac{n^{\delta}}{h(n)} \lambda \, ds\, \nu_n(du)}
\end{eqnarray}
by \cite[Lemma 5]{kaj}. By integration by parts and in view of
Potter bounds as before, for $\epsilon>0$ there exists $n_0 \in
\mathbb{N}$ such that the part of the integral for $u\geq n_0/n$ on
the right hand side of (\ref{upper}) is bounded from above by
\begin{equation}  \label{up1}
2(C+\epsilon) \int_{n_0/n}^\infty \!\! \int_{-\infty}^t \!\!
\left|\frac{\partial}{\partial u} \left[u^{1+\kappa}\!\!
\left[f\left(\frac{t-s}{u}\right)-
f\left(\frac{-s}{u}\right)\right]^{1+\kappa}\right]\right|  \lambda
\, u^{-\delta} \max\{u^{-\epsilon},u^{\epsilon}\}\, ds\, du
\end{equation}
We have
\begin{eqnarray*}
 \lefteqn{\frac{\partial}{\partial u} \left[u^{1+\kappa}\!\!
\left[f\left(\frac{t-s}{u}\right)-
f\left(\frac{-s}{u}\right)\right]^{1+\kappa}\right]  = }\\
& \displaystyle{(1+\kappa)u^{\kappa}
\left[f\left(\frac{t-s}{u}\right)-
f\left(\frac{-s}{u}\right)\right]^{1+\kappa} }\\
& \displaystyle{ + (1+\kappa)u^{\kappa}
\left[f\left(\frac{t-s}{u}\right)-
f\left(\frac{-s}{u}\right)\right]^{\kappa}
\left[-f'\left(\frac{t-s}{u}\right) \: \frac{t-s}{u}-
f'\left(\frac{-s}{u}\right)\, \frac{-s}{u}\right]}
\end{eqnarray*}
an upper bound for the absolute value of which is given by
\begin{eqnarray}  M^{1+\kappa} (1+\kappa)\left\{
[u^{-1}(u+s)^{1+\kappa}+ u^{-1}s(u+s)^{\kappa} ] 1_{R_1} +  u^\kappa
1_{R_2}   \label{up2} \right. \\  \left. + 2  u^{-1}(t-s)^{1+\kappa}
1_{R_3} +  u^{-1}t^{1+\kappa}(2+|s|)1_{R_4} \right\}\nonumber
\end{eqnarray}
by Lipschitz assumptions on $f$ and $f'$, where $R_1,\ldots , R_4$
are as in i) through iv) above, with $t_k\equiv t$. Substituting
(\ref{up2}) in (\ref{up1}) and starting the lower limit for $u$ from
0, we have an upper bound for the integral in (\ref{up1}) given by
\begin{eqnarray} \label{ints}
 \lefteqn{\int_{-\infty}^{0}\int_{-s}^{t_k -s}[(u+s)^{1+\kappa}
+s(u+s)^{\kappa}]u^{-\delta-1} \max(u^{-\epsilon},u^\epsilon)\,du \,
ds  +} \\ & \nonumber
 \displaystyle{ \int_0^{t}\int_0^{t-s} u^{\kappa-\delta}
\max(u^{-\epsilon},u^\epsilon)\, du\,  ds
 + 2\int_0^{t_k}\int_{t -s}^\infty (t-s)^{1+\kappa}
u^{-\delta-1} \max(u^{-\epsilon},u^\epsilon)\, du\,  ds }\\
& \nonumber \displaystyle{ + \int_{-\infty}^{0}\int_{t -s}^\infty
t^{1+\kappa}(2+|s|)\, u^{-\delta} \max(u^{-\epsilon},u^\epsilon)\,du
\, ds}
\end{eqnarray}
 which is finite when we choose
$\epsilon$ as in (\ref{epkappa}). On the other hand, for
$0<u<n_0/n$, we use Markov's inequality as before to get
\[
\p\{U>nu\}\, \frac{n^{\delta}}{h(n)}   \leq \mathbb{E} U\,
n_0^{\delta+\epsilon -1} u^{-\delta -\epsilon}
\]
 Then, the
finiteness of the integrals in (\ref{ints}) is sufficient again for
the integrability of a dominating function for $0<u<n_0/n<1$ which
complements (\ref{up1}). It follows from dominated convergence
theorem that the limit of  the right hand side of (\ref{upper})
exists. Therefore, possibly for $n\geq n_1$ for some $n_1\in
\mathbb{N}$, the upper bound in (\ref{upper}) is further bounded by
a multiple of its limit given by
\begin{equation} \label{t3}
C_1 \, \mathbb{E}|R|^{1+\kappa} \int_0^\infty \!\! \int_{-\infty}^t
\!\! u^{1+\kappa}\!\! \left[f\left(\frac{t-s}{u}\right)-
f\left(\frac{-s}{u}\right)\right]^{1+\kappa} \!\! \lambda
\,u^{-\delta-1} \, du\,  ds
\end{equation}
for some $C_1>2$. In view of the proof of Lemma \ref{lemma1}, the
integral in (\ref{t3}) is bounded by a constant multiple of
$t^{2+\kappa-\delta}$ which  clearly dominates
$\mathbb{E}|Z_n(t)-\mathbb{E}Z_n(t)|^{1+\kappa}$ in (\ref{upper})
for sufficiently large $n$. Since the  increments of
$\{Z_n(t)-\mathbb{E}Z_n(t):t\geq 0\}$ are stationary, this implies
that
\begin{eqnarray} \nonumber
\mathbb{E}\, [|Z_n(t_2)-\mathbb{E}Z_n(t)|^{\frac{1+\kappa}{2}}\,
|Z_n(t)-\mathbb{E}Z_n(t_1)|^{\frac{1+\kappa}{2}}] & \leq & C_2 \,
(t_2-t)^{\frac{2+\kappa -\delta}{2}}(t-t_1)^{\frac{2+\kappa -\delta}{2}} \\
& \leq & C_2 \,  (t_2 -t_1)^{2+\kappa -\delta} \label{tight}
\end{eqnarray}
for $0<t_1< t<t_2$ and some $C_2>0$, by  Cauchy-Schwarz inequality
and the assumption that $\delta<1+\kappa$. This concludes the proof
by \cite[Thm.13.5 and Eqn.(13.14)]{bill} as $2+\kappa-\delta>1$.
\finish

\begin{rmk}
  We consider the scaled measure $\nu_n$ as a scaling of the
  parameters of $\nu$. For instance, if $\nu$ was the Pareto
  distribution $\nu(du)=\delta b^\delta u^{-\delta-1} \,du$ for $u>b$, with parameters
  $\delta>0$ and $b>0$,  we
  would have $\nu_n(du)=\delta (b/n)^\delta u^{-\delta-1} \,du$ for
  $u>b/n$, which would amount to scaling the scale parameter $b$ as
  $b/n$. This
  leads to the limiting infinite measure on $\mathbb{R}_+$ as the cutoff
  parameter $b$ decreases. Although similar interpretations are possible for other
  probability measures $\nu$ as well, the scaling
  $\nu_n(du)=\nu(d(nu))$ is essentially a time scaling. The random
  variable $U$, the duration of the sessions, has time
  interpretation. We look at the workload of individual pulses
  over shorter time periods by scaling $U$ as $U/n$.
\end{rmk}

Note that the limiting process given in Theorem \ref{ICR},
specifically with $f$ of (\ref{increasingf}), is obtained as a limit
of sums of scaled renewal processes in \cite{gaigalas} as shown in
\cite{gaig}.

\section{Fractional Brownian Motion Limit}

 Fractional
Brownian motion   is a mean zero Gaussian process $Z$ on
$\mathbb{R}_+$ with $Z(0)=0$ and covariance
\[
\mbox{Cov}(Z(t_1),Z(t_2))=\frac{\sigma^2}{2}\,
(|t_1|^{2H}+|t_2|^{2H}-|t_1-t_2|^{2H}) \qquad t_1,t_2\geq 0
\]
for $t_1,t_2\geq 0$, $\sigma>0$ and Hurst parameter $0<H<1$
\cite{samor}.

In this section, we scale the log-price process as follows to
approximate a fractional Brownian motion in the limit. Let the rate
$R$ be scaled as $R/n$ which can be interpreted as a decrease in the
effect of an order in absolute value as $n$ increases. This decrease
may have arised from an underlying decrease in the volume of the
transaction, for example. On the other hand, we will let the arrival
rate $\lambda$ of orders increase with a factor which is a function
of $n$. For each $n\in\mathbb{Z}_+$, let
 $N_n$ denote the Poisson random measure with scaled mean measure $\mu_n$ that involves the
 scaled  arrival rate and possibly further scalings.
 In the following theorems, we prove convergence of
 $Z_n -\mathbb{E}Z_n$ to fBm
 with a properly scaled measure $\nu_n$  for
  a finite measure $\nu$ as in (\ref{nu})
   and  with $\nu$ as in (\ref{mandelu}). Compensated Poisson random
measure is used when $\mathbb{E}Z_n$ does not exist.

\begin{lem} Let $(x_n)$ be a strictly positive real sequence and
$c\in\mathbb{R}$. If $\lim(x_n)=\infty$, then
\[
\lim_{n\rightarrow \infty} x_n^2 (e^{i c
/x_n}-1-ic/x_n)=-\frac{c^2}{2}
 \; .
\]  \label{lemma 2}
\end{lem}
\prf Consider the Taylor expansion \cite[pg.184-186]{complex} of the
function $h(x)=e^{i c /x}$ around 0  on the disk
$D=\{x\in\mathbb{R}:|x|<2\}$. We have
\[
x_n^2 \left| (e^{i c /x_n}-1-ic/x_n)+\frac{c^2}{2x_n^2}\right| =
x_n^2 |R_2(x_n)|
\]
where $R_2 $ is the remainder after 3 terms. For fixed $r$ with
$1<r<2$, we can express $R_2 $ as
\[
R_2(x_n)=\frac{1}{2\pi
i\,x_n^3}\int_{|x|=r}\frac{h(x)}{x^3(x-1/x_n)}\,dx
\]
By  continuity of $h$, there exists $M>0$ such that $|h(x)|\leq M$
for all $x\in D$, in particular for $|x|=r$.  Clearly,
$|x-1/x_n|\geq |x|-|1/x_n|=r-1/x_n$ since $x_n$ is positive. Then,
we have the following bound for $|R_2(x_n)|$:
\[
|R_2(x_n)|\leq \frac{M\, r}{r-1/x_n}\left(\frac{1/x_n}{r}\right)^3.
\]
Taking the limit as $n\rightarrow\infty$, we obtain
$x_n^2|R_2(x_n)|\rightarrow 0$.
 \finish

\begin{thm}  \label{FCR} Suppose that  $\nu$ is a probability measure
 satisfying $\int u \, \nu(du)<\infty$ and has a regularly varying
 tail as given in
 (\ref{nu}), the function $f:\mathbb{R}\rightarrow \mathbb{R}$ is Lipschitz
continuous with $f(x)=0$ for all $x\leq 0$, $f(x)=f(1)$ for all
$x\geq 1$ and is also differentiable with $f'$ satisfying  a
Lipschitz condition a.e., and   $\mathbb{E}R^2<\infty$. Let
\[
Z_n(t)=\int_{-\infty}^{\infty}\int_0^{\infty}
\int_{-\infty}^{\infty}
\frac{r}{n}\,u\left[f\left(\frac{t-s}{u}\right)-
f\left(\frac{-s}{u}\right)\right]N_n(ds,du,dr)
\]
and
\[ \mu_n (ds,du,dr)= \frac{n^{2+\delta}}{h(n)}\lambda \, ds \, \nu_n(du)\,
\gamma(dr)
\]
where $\nu_n(du) = \nu(d(nu))$ and $1<\delta<2$. Then, the process
$\{Z_n(t)-\mathbb{E} Z_n(t) ,t\geq 0\}$  converges in law to an fBm
with variance parameter
$$\sigma^2=\lambda\, \mathbb{E}R^2\int_{-\infty}^{\infty}\int_0^{\infty} \left[
f\left(\frac{1-s}{u}\right)- f\left(\frac{-s}{u}\right)\right]^2
u^{1-\delta} du \, ds $$ as $n\rightarrow\infty$.
\end{thm}

 \prf The characteristic function
  $\mathbb{E}\exp i\sum_{k=1}^m\xi_k    [Z_n(t_k)-\mathbb{E}
Z_n(t_k)]  $ for $\xi_k\in \mathbb{R}$, $t_k\geq 0$ and $m\in
\mathbb{N}$ is given by
\begin{eqnarray}
&& \exp \int_{-\infty}^{\infty}\int_0^{\infty}
\int_{-\infty}^{\infty}\Bigg\{
e^{i\sum_{k=1}^m\xi_k\frac{r}{n}u\left[f\left(\frac{t_k-s}{u}\right)-
f\left(\frac{-s}{u}\right)\right]}-1\Bigg. \nonumber\\
&&\qquad
\left.-i\sum_{k=1}^m\xi_k\,\frac{r}{n}\,u\left[f\left(\frac{t_k-s}{u}\right)-
f\left(\frac{-s}{u}\right)\right]\right\}
\frac{n^{2+\delta}}{h(n)}\: \lambda\,
 ds\, \nu_n( du)\,\gamma(dr)\label{chr}.
\end{eqnarray}
The same approach will be followed as in the proof of Theorem
\ref{ICR}. By integration by parts, we find that the exponent of
(\ref{chr}) is given by
\begin{equation}
\int \int \int  \partial_u g(s,u,r/n) \, \p\{U>nu\}\,
\frac{n^{2+\delta}}{h(n)}\: \lambda \, ds \, du\, \gamma(dr)\: .
\label{byparts2}
\end{equation}
Using  Potter bounds \cite{bingham} and Lipschitz conditions on $f$
and $f'$, we get an inequality similar to (\ref{bound}) for $u\geq
n_0/n$ given by
\begin{eqnarray} \label{boundFCRa}
\lefteqn{\int  \int_{n_0/n}^\infty \int |\partial_u g(s,u,r/n)| \,
\p\{U>nu\}\,
\frac{n^{2+\delta}}{h(n)} \, \lambda \, ds \, du\, \gamma(dr)} \\
&  \leq  4 M^2 (C+\epsilon) \displaystyle{\sum_k |\xi_k|
 \!\! \int \!\! \int_{0}^\infty  \!\!\!\! \int
|r| \, \tilde{B}(s,u,t_k) \max (u^{-\epsilon},u^{\epsilon})\,
 u^{-\delta}\, \lambda  \, ds \, du \, \gamma(dr) } \nonumber
\end{eqnarray}
where $\tilde{B}$ is similar to (\ref{B}) but with $\kappa=1$ by
hypothesis, and $\epsilon >0$ and $n_0 \in \mathbb{N}$. Precisely,
\begin{eqnarray*}
\lefteqn{\tilde{B}(s,u,t_k)= \left( 1\wedge \sum_j    M |\xi_j|
\, |r| \,u \, 1_{\{s\leq t_j\}} \right)} \nonumber \\
 & \qquad   . \: \left[  1_{R_{1,k}}+   1_{R_{2,k}}+ 2 \frac{t_k-s}{u}
 \, 1_{R_{3,k}}+  \left(2   \frac{t_k}{u}   +
\frac{|s| t_k}{u^2} \right) 1_{R_{4,k}} \right]
\end{eqnarray*}
 If we choose $\epsilon>0$ such that
\[
1<\delta -\epsilon<\delta <\delta+\epsilon<2 \: ,
\]
then the right hand side of (\ref{boundFCRa}) is  finite  along the
same lines   of the proof of Theorem \ref{ICR} with $\kappa=1$. On
the other hand, we can bound (\ref{byparts2}) for $0<u\leq 1$
similarly. Therefore, we can use dominated convergence theorem. We
have
\[
\lim_{n\rightarrow \infty} \p\{U>nu\}\, \frac{n^{\delta}}{h(n)}
=\frac{u^{-\delta}}{\delta}
\]
as in (\ref{kucuklimit}), and
\[
\lim_{n\rightarrow \infty} n^2\partial_u g(s,u,r/n)  =  \partial_u
\lim_{n\rightarrow \infty}n^2 g(s,u,r/n)
\]
as $g$ is bounded, hence, uniformly continuous. Then, we get
\begin{eqnarray}
\lefteqn{\lim_{n\rightarrow \infty}n^2 g(s,u,r/n) =} \label{limg} \\
&\displaystyle{-\frac{1}{2}\sum_{k=1}^m\sum_{j=1}^m\xi_j\xi_k r^2
\left[f\left(\frac{t_j-s}{u}\right)
-f\left(\frac{-s}{u}\right)\right] \nonumber
\left[f\left(\frac{t_k-s}{u}\right)-f\left(\frac{-s}{u}\right)\right]u^2}
\end{eqnarray}
 by Lemma \ref{lemma 2}. We now
  revert (\ref{byparts2}) after the limits above, by another integration by
  parts,
  and
get the limit of (\ref{chr}) as
\begin{eqnarray*}\!\!\! \!\!\! \exp&& \!\!\! \!\!\! \!\!\!
\left\{  -\int_{-\infty}^\infty
\right.\int_0^{\infty}\int_{-\infty}^{\infty}
\frac{1}{2}\sum_{k=1}^m\sum_{j=1}^m
\xi_j\xi_k r^2 \nonumber\\
&& \!\!\! \!\!\! \left.\left[f\left(\frac{t_j-s}{u}\right)
-f\left(\frac{-s}{u}\right)\right]
\left[f\left(\frac{t_k-s}{u}\right)-f\left(\frac{-s}{u}\right)\right]u^{2}
\lambda \, ds \, u^{-\delta-1}\, du \, \gamma(dr)\right\}
\end{eqnarray*}
which is the characteristic function of $\sum_{k=1}^m\xi_kZ(t_k)$,
where $Z=(Z(t_1),\ldots,Z(t_m))$ is a Gaussian vector with zero mean
and covariance
\begin{equation}
\lambda\mathbb{E}R^2\int_0^{\infty}\int_{-\infty}^{\infty}\left[f\left(\frac{t_j-s}{u}\right)
-f\left(\frac{-s}{u}\right)\right]
\left[f\left(\frac{t_k-s}{u}\right)-f\left(\frac{-s}{u}\right)\right]
u^{2}ds \,u^{-\delta-1}\, du. \label{covariance_matrix}
\end{equation}
When (\ref{covariance_matrix}) is evaluated at
 $t_j=t_k=1$,
 the variance coefficient is found
to be $$ \sigma^2: =
\mbox{Var}(Z(1))=\lambda\mathbb{E}R^2\int_0^\infty\int_{-\infty}^\infty\left[f\left(\frac{1-s}{u}\right)-
f\left(\frac{-s}{u}\right)\right]^2 \,ds\,u^{1-\delta}\, du
$$
which is finite by Lemma \ref{lemma1} with $\kappa=1$. Using the
identity $2ab= -(a-b)^2+a^2+b^2$ for $a,b\in \mathbb{R}$ and making
several change of variables, we  find that the covariance of $Z$  in
(\ref{covariance_matrix}) is given by
$$\mbox{Cov}(Z(t_j),Z(t_k))=\frac{\sigma^2}{2}\, (t_j^{2H}+t_k^{2H}-|t_j-t_k|^{2H})$$
for $t_j,t_k\geq 0$ with $H=(3-\delta)/2$. By definition, $Z$ has
the characteristic function of an fBm.

   Convergence in the Skorohod topology on $D(0,\infty)$ follows
   along the same lines of proof of Theorem \ref{ICR}. In this case, we have
\[
\mathbb{E}|Z_n(t)-\mathbb{E}Z_n(t)|^2 = \mathbb{E}R^2 \int_0^\infty
\!\! \int_{-\infty}^t \!\! u^2\!\!
\left[f\left(\frac{t-s}{u}\right)-
f\left(\frac{-s}{u}\right)\right]^2 \!\! \lambda \, ds\, \nu_n(du)
\]
and (\ref{tight}) holds with   $\kappa=1$. \finish

 As an example, the continuous flow rate model studied in
\cite{kaj} is given by
\begin{equation}
Z(t)=\int_{-\infty}^{\infty}\int_{0}^\infty\int_{-\infty}^{\infty}[(t-s)^+\wedge
u-(-s)^+\wedge u]\,r \, N(ds,du,dr)
 \label{taqqu5}
 \end{equation}
with $f$ replaced by the special form (\ref{increasingf}) in
(\ref{Zf}). In \cite[Thm.2]{kaj}, the limit is  studied when the
speed of time increases in proportion to the intensity of Poisson
arrivals. To balance the increasing trading intensity $\lambda_n$,
time is speeded up by a factor $n$ and  the size is normalized by
a factor $\lambda_n^{1/2}n^{(3-\delta)/2}$ provided that
$\lambda_n/n^{\delta -1}\rightarrow \infty$. We can let
$\lambda_n= n^{\varepsilon+\delta-1}$ with $\varepsilon>0$. Taking
$\varepsilon=2$, we show the equivalence of the scaling of
\cite[Thm.2]{kaj} to the scaling in Theorem \ref{FCR}. Note that
$\lambda_n= n^{1+\delta}$. The scaled and centered process has the
form
\begin{eqnarray}
  \lefteqn{\frac{Z(nt)-\mathbb{E}Z(nt)}{\lambda_n^{1/2}n^{(3-\delta)/2}}} \nonumber \\
   && =\frac{1}{ n^{ 2} }\:
\int_{-\infty}^\infty\int_0^\infty\int_{-\infty}^\infty
 ru  \left[f\left(\frac{nt-s}{u}\right)- f\left(\frac{-s}{u}\right)\right]  \,\tilde{N}_n(ds,du,dr)
 \label{taqqu22}
\end{eqnarray}
where we have written an effect function $f$ in  general. Then, we
can make change of variables $s \rightarrow ns$ and $u \rightarrow
nu$ to get
\begin{eqnarray}
\lefteqn{\frac{Z(nt)-\mathbb{E}Z(nt)}{n^{ 2} }} \nonumber \\
   && =\frac{1}{ n^{ 2} }\:
\int_{-\infty}^\infty\int_0^\infty\int_{-\infty}^\infty
 r \, nu  \left[f\left(\frac{nt-ns}{nu}\right)- f\left(\frac{-ns}{nu}\right)\right]
 \,\tilde{N}_n(d(ns),d(nu),dr) \nonumber\\
 && =
\int_{-\infty}^\infty\int_0^\infty\int_{-\infty}^\infty
 \frac{r}{n} \, u  \left[f\left(\frac{t-s}{u}\right)- f\left(\frac{-s}{u}\right)\right]
 \,\tilde{N}_n(d(ns),d(nu),dr) \label{equivform}
\end{eqnarray}
where the mean measure is
\begin{eqnarray}
\mu_n(d(ns), d(nu),dr)&=&\lambda_n\,(nds)\,\nu( d(nu)) \gamma(dr) \nonumber\\
&=& n^{2+\delta}   \, ds \,\nu( d(nu))\gamma(dr)\: .
\label{scalmean}
\end{eqnarray}
In Theorem \ref{FCR}, we start with the scaled process
(\ref{equivform}) essentially. This can be observed by the fact that
\[
 N_n(d(ns),d(nu),dr) \stackrel{d}{=} N_n'(ds,du,dr)
\]
for a Poisson random measure $N'$ with mean measure
\[
\mu_n'(ds,du,dr)=\mu_n(d(ns), d(nu),dr)
\]
  by definition of a Poisson random measure
\cite{kall},\cite[Def.V.2.2]{cinlarlec1}. Equivalence of the
scalings in Theorem \ref{FCR} and \cite[Thm.2]{kaj} is in
distributional sense. However, this is sufficient for equivalence as
the convergence results are in distribution rather than almost sure
sense.
 Therefore, we can apply Theorem
\ref{FCR} to obtain the limit as a fBm with variance parameter
\[\sigma^2= \mathbb{E}R^2 \int_0^\infty\int_{-\infty}^\infty
[(1-s)^+\wedge u-(-s)^+\wedge u]^2 \, ds\, u^{-\delta-1}du=
\frac{\mathbb{E}R^2}{(2-\delta)(3-\delta)}\: . \]

It is shown in \cite{kaj} that the asymptotic behavior of the
ratio $\lambda_n/n^{\delta-1}$  determines the type of the limit
process when time is speeded up by a factor $n$. For a choice of
sequences $\lambda_n$ and $n$, the random variable
$\sharp(\lambda_n,n)$ denotes the number of effects still active
at time n. It measures the amount of very long pulses that are
alive and how much they contribute to the total price. The
expected value of the random variable $\sharp(\lambda_n,n)$ is
\[
\mathbb{E}\,\,\sharp(\lambda_n,n)\sim\frac{\lambda_n}{n^{\delta-1}}
\]
for large $n$. The limit is considered in the cases where this value
tends to a finite positive constant, to infinity, or to zero as
$\lambda_n$ and $n$ go to infinity. We have already studied the case
of finite constant in Theorem \ref{ICR} and infinity in Theorem
\ref{FCR}, the so-called intermediate and fast connection rates,
respectively, in view  of telecommunication applications. As shown
above, our scalings do not involve time scaling. They can be
physically understood as scalings of the parameters of the log-price
process. The slow connection rate will be investigated similarly in
terms of the model parameters in Theorem \ref{SCR} in the next
section.

The next theorem is a simpler  version of Theorem \ref{FCR} due to
the form of the  measure $\nu$.
 Note that
(\ref{scalmean}) can be approximated as
\begin{eqnarray*}
 \mu_n(d(ns),d(nu), dr) & \sim  & n^{2+\delta} \,ds\,n^{-\delta}u^{-\delta-1} du\,\gamma(dr)
\\& = & n^{2}u^{-\delta-1}\,  ds\,du\,\gamma(dr)
\end{eqnarray*}
for large $n$. This scaling is used below with the simpler form of
$\nu$. It can be interpreted as half way in taking the more involved
limit of Theorem \ref{FCR}.

\renewcommand{\theenumi}{\roman{enumi}}
\begin{thm} \label{Manthm} Let
\[
\tilde{Z}_n(t)=\int_{-\infty}^{\infty}\int_0^{\infty}
\int_{-\infty}^{\infty}
\frac{r}{n}\,u\left[f\left(\frac{t-s}{u}\right)-
f\left(\frac{-s}{u}\right)\right]\tilde{N}_n(ds,du,dr)
\]
where $\tilde{N}=N-\mu$ and
\[ \mu_n(ds,du,dr)=n^2\lambda \,
u^{-\delta-1}\, ds\, du\, \gamma(dr).
\]
Suppose that $\mathbb{E}R^2<\infty$ and $f:\mathbb{R}\rightarrow
\mathbb{R}$ is a Lipschitz continuous function satisfying either of
the following conditions
\begin{enumerate}
\item  $f(x)=0$ for all $x\leq 0$   and $f(x)=f(1)$ for all $x\geq
1$,  or
\item $f$ has a compact support.
\end{enumerate}
 Then, the process $\{\tilde{Z}_n(t), t\geq 0\}$, for $1<\delta<3$,
  converge in law to an fBm with variance parameter
$$\sigma^2=\lambda\, \mathbb{E}R^2\int_{-\infty}^{\infty}\int_0^{\infty}\left[
f\left(\frac{1-s}{u}\right)- f\left(\frac{-s}{u}\right)\right]^2
u^{1-\delta} du \, ds$$ as \  $n\rightarrow\infty$.
\end{thm}
\noindent \textbf{Proof:}  Although it can be found from the
characteristic function of $Z(t)$ in Proposition \ref{prop2} that
$\mathbb{E} Z(t) =0$ for all $t\geq 0$ under assumption ii, we form
$\tilde{Z}$ as above since  $\mathbb{E} Z(t)$ may not exist with
assumption i. For the convergence of finite dimensional
distributions of $\{\tilde{Z}_n(t), t\geq 0\}$, consider the
characteristic function $\mathbb{E}\exp i\sum_{k=1}^m\xi_k
\tilde{Z}_n(t_k) $ for $\xi_k\in \mathbb{R}$, $t_k\geq 0$ and $m\in
\mathbb{N}$. It is given by
\begin{equation}
  \exp \int_{-\infty}^{\infty}\int_0^{\infty}
\int_{-\infty}^{\infty} g(s,u,r/n) n^2\lambda\,u^{-\delta-1}\,
 ds\, du\,\gamma(dr) \label{mandelbrot2.5}
\end{equation}
where $g$ is given in (\ref{g}).  Note that the characteristic
function exists since $\tilde{Z}_n(t_k)$ are well defined in view of
(\ref{secmoment}) which follows from Lemma \ref{lemma1} with
$\kappa=1$ under assumption i, and by Proposition \ref{prop2} under
assumption ii. As $n\rightarrow\infty$, we will show that the above
characteristic function converges to
\begin{eqnarray}\!\!\! \!\!\! \exp&& \!\!\! \!\!\! \!\!\!
\left\{  -\int_{-\infty}^\infty
\right.\int_0^{\infty}\int_{-\infty}^{\infty}
\frac{1}{2}\sum_{k=1}^m\sum_{j=1}^m
\xi_j\xi_k r^2 \nonumber\\
&& \!\!\! \!\!\! \left.\left[f\left(\frac{t_j-s}{u}\right)
-f\left(\frac{-s}{u}\right)\right]
\left[f\left(\frac{t_k-s}{u}\right)-f\left(\frac{-s}{u}\right)\right]u^{2}
 \lambda  \,u^{-\delta-1} \, ds \,du \, \gamma(dr)\right\} \label{mandelbrot_bound}
\end{eqnarray}
 Due to the
inequality $|e^{ix}-1-ix|<\frac{1}{2}x^2$ for $x\in\mathbb{R}$, the
integrand in (\ref{mandelbrot_bound})  is an upper bound to
$$|g(s,u,r/n)|\, n^2$$
 Therefore,
dominated convergence theorem allows us to take the limit inside the
integral in (\ref{mandelbrot2.5}). That is, we must find
\[\lim_{n\rightarrow\infty}
\left(e^{i\sum_{k=1}^m\xi_k\frac{r}{n}u\left[f\left(\frac{t_k-s}{u}\right)-
f\left(\frac{-s}{u}\right)\right]}-1
-i\sum_{k=1}^m\xi_k\frac{r}{n}u\left[f \left(\!\frac{t_k-s}{u} \!
\right)  -  f   \left( \! \frac{-s}{u}\! \right) \right]\right)n^2
\]
which is now equal to
\begin{equation}-\frac{1}{2}\sum_{k=1}^m\sum_{j=1}^m\xi_j\xi_k r^2
\left[f\left(\frac{t_j-s}{u}\right)
-f\left(\frac{-s}{u}\right)\right]
\left[f\left(\frac{t_k-s}{u}\right)-f\left(\frac{-s}{u}\right)\right]u^2
\label{limit}
\end{equation}
by Lemma \ref{lemma 2}. This shows that (\ref{mandelbrot2.5})
converges to (\ref{mandelbrot_bound})  as $n\rightarrow\infty$ by
the continuity of the exponential function. The variance is
evaluated in the proof of Theorem \ref{FCR}.

   To complete the proof, we need to show convergence in
   $D(0,\infty)$ with Skorohod topology. This is straight forward
    since the variance of $\tilde{Z}_n(t)$ is
    already free of $n$ and is bounded by a constant multiple of $t^{3-\delta}$
    by the proof of Lemma \ref{lemma1}.
 \finish

\begin{rmk}
Theorem \ref{Manthm} with condition ii. is Theorem 3.1 of
\cite{mandel} where it is noted that a fractional Brownian motion
with $H>1/2$ can  be approximated if the pulse is continuous and has
a compact support. Condition i. above considers an effect function
which is continuous, but with no compact support as an alternative.
\end{rmk}


\section{L\'{e}vy Process Limit}

   A process with stationary and independent increments is called
a L\'{e}vy process \cite{bert,sato}. The results of this section
concerns a particular class of L\'{e}vy processes, namely stable
L\'{e}vy motion \cite{samor}.  Let
$\mathbb{R}_0=\mathbb{R}\setminus\{0\}$, and let $\delta \in (1,2)$,
 and $\beta \in [-1,1]$ be  the index of
stability and skewness parameter, respectively. Then, a
$\delta$-stable L\'{e}vy motion $L$ with mean 0 can be defined
through its characteristic function
\[
\mathbb{E}e^{i\xi L(t)} =
\exp\{-t\,\sigma^\delta\,|\xi|^\delta[1-i\beta(\mbox{sign}\, \xi)
\tan (\pi\delta/2)]\}
\]
for $\xi \in \mathbb{R}$, where $\sigma\geq 0$ is a scale parameter.


 In this section, we prove that the limiting process is a $\delta$-stable L\'{e}vy
motion under different scalings of the price process. Theorem
\ref{SCR} considers a probability measure $\nu$ and Theorem
\ref{SCRMandelbrot} starts with its limiting form. For simplicity of
notation, we take $f(1)=1$ in the effect function.

\begin{lem} \label{lemma3} Let $N$ be a Poisson random measure with mean measure
\[
\mu=\lambda \, ds \, u^{-\delta -1}du \,\gamma(dr) \]
 and
$\tilde{N}=N-\mu$. Then,
\[
\int_{-\infty}^{\infty}\int_0^{\infty} \int_{0}^{t} r\, u \,
\tilde{N}(ds,du,dr) \stackrel{d}{=} (\lambda \,C_1)^{1/\delta}\,
L_1(t) + (\lambda \,C_2)^{1/\delta} \, L_2(t)
\]
where $C_1=\int_0^{\infty} r^{\delta}\gamma(dr)$,
$C_2=\int_{-\infty}^0 |r|^{\delta}\gamma(dr)$, and $L_1$ and $L_2$
are  independent $\delta$-stable L\'{e}vy motions with mean 0,
skewness intensity $\beta$ equal to 1 and $-1$, respectively, and
scale parameter
\[
\sigma=\left[-\frac{2\Gamma(2-\delta)}{\delta(\delta-1)}\cos
\frac{\pi\delta}{2}\right]^{1/\delta} \; .
\]
\end{lem}

\prf Putting $u'=|r|u$, we get
\begin{eqnarray}
\lefteqn{\int_{-\infty}^{\infty}\int_0^{\infty} \int_{0}^{t} r\, u
\, \tilde{N}(ds,du,dr)} \nonumber \\ & \displaystyle{=
\int_0^{\infty} \int_0^{\infty} \int_{0}^{t}   u' \,
[N(ds,d(u'/r),dr)-\lambda \, ds \, u'^{-\delta -1}r^{\delta}
\, du' \,\gamma(dr)]} \nonumber \\
& \displaystyle{ \qquad
 - \int_{-\infty}^{0}\int_0^{\infty}
\int_{0}^{t} u' \, [N(ds,d(u'/|r|),dr)-\lambda \, ds \,
u'^{-\delta -1} |r|^{\delta}\, du' \,\gamma(dr)]}  \nonumber \\
& \!\!\!\!\!\! \!\!\!\!\!\! \!\!\!\!\!\! \!\!\!\!\!\! \!\!\!\!\!\!
\!\!\!\!\!\!  \!\!\!\!\!\! \!\!\! \displaystyle{= \int_0^{\infty}
\int_{0}^{t} u' \label{star}
\, [N_1'(ds,du')-C_1 \lambda \, ds \, u'^{-\delta -1}\, du'  ]}\\
 & \!\!\!\!\!\! \!\!\!\!\!\! \!\!\!\!\!\! \!\!\!\!\!\! \!\!\!\!\!\!
\!\!\!\!\!\!  \!\!\!\!\!  \displaystyle{  - \int_0^{\infty}
\int_{0}^{t} u' \, [N_2'(ds,du')-C_2 \lambda \, ds \, u'^{-\delta
-1} \, du' ]} \nonumber
\end{eqnarray}
where we have put
\begin{eqnarray*}
N_1' (ds,du')& := & \int_{\{r:r \in (0,\infty)\}}  N(ds,d(u'/r),dr)              \\
N_2' (ds,du')& := & \int_{\{r:r\in(-\infty,0)\}}N(ds,d(u/|r|),dr)
\end{eqnarray*}
and
\[ C_1 := \int_0 ^\infty r^{\delta}\gamma(dr), \;\;\;
C_2 := \int_{-\infty}^0 |r|^{\delta}\gamma(dr),
\]
 It is easy to verify that both $N_1'$ and $N_2'$ are Poisson
random measures with means $C_i\lambda ds \, u'^{-\delta -1} \,
du'$, $i=1,2$, respectively. What is more, they are independent as
their domains are disjoint. Making another change of variable $u=
u'/(\lambda C_i)^{1/\delta}$ for $i=1,2$ in respective integrals in
(\ref{star}), we get
\begin{eqnarray}
\lefteqn{\int_{-\infty}^{\infty}\int_0^{\infty} \int_{0}^{t} r\, u
\, \tilde{N}(ds,du,dr)} \nonumber \\
& = \displaystyle{(\lambda C_1)^{1/\delta}  \int_0^{\infty}
\int_{0}^{t} u\, \nonumber
[N_1'(ds,d((\lambda C_1)^{1/\delta}u))-  ds \, u^{-\delta -1}\, du]} \\
& \qquad \displaystyle{-(\lambda C_2)^{1/\delta} \int_0^{\infty}
\int_{0}^{t} u\, [N_2'(ds,d((\lambda C_2)^{1/\delta}u))- ds \,
u^{-\delta
-1}\, du]} \nonumber \\
&  \!\!\!\!\!\! \!\!\!\!\!\! \!\!\!\!\!\! \!\!\!\!\!\! \!\!\!\!\!\!
 \!\! \!\!\!\!\!\!\!\!\!\!\!\! \!\!\!\!\!\!  \!\!\!\!  \!\!\!\! \!\!
 =:\displaystyle{(\lambda C_1)^{1/\delta}    L_1'(t) -(\lambda
C_2)^{1/\delta}   L_2'(t)}
\end{eqnarray}
where $L_1'$ and $L_2'$ are independent $\delta$-stable L\'{e}vy
motions with  skewness parameter $\beta=1$ and scale parameter
\[
\sigma=\left[-\frac{2\Gamma(2-\delta)}{\delta(\delta-1)}\cos
\frac{\pi\delta}{2}\right]^{1/\delta}
\]
for $\xi \in \mathbb{R}$ \cite[pg.s 5,156,350]{samor}. Now, we have
\[
-L_2'\stackrel{d}{=}L_2:=\int_{-\infty}^0 \int_{0}^{t} u\,
[N_2'(ds,d(-(\lambda C_2)^{1/\delta}u))- ds \, |u|^{-\delta -1}\,
du]
\]
where $L_2$ is also a $\delta$-stable L\'{e}vy motion since
$N_2''(ds,du):=N_2'(ds,d(-(\lambda C_2)^{1/\delta}u))$ is a Poisson
random measure on $\mathbb{R}_{+}\times \mathbb{R}_-$ with mean
measure $ds \, |u|^{-\delta -1}\, du$, but skewness parameter
$\beta=-1$ \cite[pg.5]{samor}. We can take $L_1=L_1'$ and the result
follows.
 \finish

\begin{thm} \label{SCR}  Suppose that  $\nu$ is a probability measure
 satisfying $\int u \, \nu(du)<\infty$ and has a regularly varying
 tail as given in
 (\ref{nu}), the function $f:\mathbb{R}\rightarrow \mathbb{R}$ is Lipschitz
continuous with $f(x)=0$ for all $x\leq 0$, $f(x)=f(1)=1$ for all
$x\geq 1$ and is also differentiable with $f'$ satisfying  a
Lipschitz condition a.e., and $\mathbb{E}|R|^{1+\kappa}<\infty$
for some $0<\kappa\leq 1$
such that $1+\kappa>\delta>1$. Let
\[
Z_n(t)=\int_{-\infty}^{\infty}\int_0^{\infty}
\int_{-\infty}^{\infty} n^{1-\alpha/\delta}r
\,u\left[f\left(\frac{t-s}{u}\right)-
f\left(\frac{-s}{u}\right)\right]N_n(ds,du,dr)
\]
and
\[ \mu_n (ds,du,dr)= \frac{n^{\alpha}}{h(n^{\alpha/\delta})}\lambda \, ds \, \nu_n(du)\,
\gamma(dr)
\]
where $\nu_n(du) = \nu(d(nu))$ and  $0<\alpha<\delta$.
Then, the process $\{Z_n(t)-\mathbb{E} Z_n(t), t\geq 0\}$, for
$1<\delta<2$, converges in law to
\[
(\lambda \,\mathbb{E} R^{\delta}1_{\{R>0\}})^{1/\delta}\, L_1(t) +
(\lambda \,\mathbb{E} |R|^{\delta}1_{\{R<0\}})^{1/\delta} \, L_2(t)
\]
 as $n\rightarrow\infty$, where $L_1$ and $L_2$
are  independent $\delta$-stable L\'{e}vy motions with mean 0, and
skewness intensity $1$ and $-1$, respectively.
\end{thm}

\prf The characteristic function of the scaled process is given by
\[
  \exp \int_{-\infty}^{\infty}\int_0^{\infty}
\int_{-\infty}^{\infty} g(s,u,n^{1-\alpha/\delta}r )
\frac{n^{\alpha}}{h(n^{\alpha/\delta})}\: \lambda\,
 ds\, \nu_n( du)\,\gamma(dr)
\]
where $g$ is as  in (\ref{g}). By integration by parts, we get
\begin{equation}
\exp \int_{-\infty}^{\infty}\int_0^{\infty} \int_{-\infty}^{\infty}
\partial_u g(s,u,n^{1-\alpha/\delta}r )\p\{U>nu\}
\frac{n^{\alpha}}{h(n^{\alpha/\delta})}\: \lambda\,
 ds\,  du \,\gamma(dr)\label{chr2}
\end{equation}
since $\nu_n(du) = \nu(d(nu))$. Making a change of variable $u$ to
$u/n^{1-\alpha/\delta}$, the logarithm of (\ref{chr2}) is equal to
\begin{equation} \label{byparts3} \int \int
\int \frac{1}{n^{1-\alpha/\delta}}\, \partial_u
g(s,u/n^{1-\alpha/\delta},n^{1-\alpha/\delta}r) \,
\p\{U>n^{\alpha/\delta}u\}\,
\frac{n^{\alpha}}{h(n^{\alpha/\delta})}\: \lambda \, ds \, du \,
\gamma(dr)\: .
\end{equation}
Using  Potter bounds  and Lipschitz conditions on $f$ and $f'$, we
get an inequality similar to (\ref{bound}). We can bound
$\p\{U>n^{\alpha/\delta}u\}\,  n^{\alpha}/h(n^{\alpha/\delta})$ as
in (\ref{14}) and consider $|\partial_u
g(s,u/n^{1-\alpha/\delta},n^{1-\alpha/\delta}r)|/n^{1-\alpha/\delta}$
separately. As a result, for fixed $\epsilon>0$, there exists
$n_0\in \mathbb{N}$ such that for all $n$ with
$n^{\alpha/\delta}\geq n_0$,   we have the following upper bound for
the absolute value of (\ref{byparts3}) when evaluated over $u\geq
n_0/n^{\alpha/\delta}$
\begin{eqnarray} \label{boundSCR}
4 M  (C+\epsilon) \displaystyle{\sum_k |\xi_k|
 \!\! \int \!\! \int_{0}^\infty  \!\!\!\! \int
|r| \, B'(s,u,t_k,n) \max (u^{-\epsilon},u^{\epsilon})\,
 u^{-\delta}\, \lambda  \, ds \, du \, \gamma(dr) } \nonumber
\end{eqnarray}
where   ${B}'$ is  analogous  to (\ref{B}) satisfying
\begin{eqnarray} \label{42}
B'(s,u,t_k,n)\leq   1_{R_{1,k,n}}+\left( 1\wedge \sum_j M
|\xi_j|^\kappa |r|^\kappa u^\kappa 1_{\{s\leq t_j \}}\right).\,
1_{R_{2,k,n}} \\ \nonumber +2\,
\frac{t_k-s}{u/n^{1-\alpha/\delta}}\: 1_{R_{3,k,n}}+
\left(2\frac{t_k}{u/n^{1-\alpha/\delta}} +
\frac{|s|t_k}{u^2/n^{2(1-\alpha/\delta)}}\right)1_{R_{4,k,n}}
\end{eqnarray}
and  $R_{1,k,n}\ldots, R_{4,k,n}$ are analogous to $R_{1,k}, \ldots,
R_{4,k}$  with $u$ replaced by $u/n^{1-\alpha/\delta}$. The right
hand side of (\ref{boundSCR}) is integrable when $\epsilon$ is
chosen as in (\ref{epkappa}) as shown in Appendix B.

For $u< n_0/n^{\alpha/\delta}$, we can find a dominating function
for the integrand in (\ref{byparts3}) using Markov's inequality. As
in the proof of Theorem \ref{ICR}, we have
\begin{eqnarray} \nonumber
\lefteqn{\frac{1}{n^{1-\alpha/\delta}}|\partial_u
g(s,u/n^{1-\alpha/\delta},n^{1-\alpha/\delta}r)| \,
\p\{U>n^{\alpha/\delta}u\}\,
\frac{n^{\alpha}}{h(n^{\alpha/\delta})}}\\ &
\;\;\;\;\;\;\;\;\;\;\;\;\;\;\;\; \displaystyle{ \leq 2M \frac{\E
U}{u} \frac{n^{\alpha-\alpha/\delta}}{h(n^{\alpha/\delta})} \, |r|
\sum_k|\xi_k| B'(s,u,t_k,n) } \label{newbound}
\end{eqnarray}
where $B'$ satisfies (\ref{42}). Using (\ref{42})  and using an
inequality similar to (\ref{bound3}) in view of the assumption $u<
n_0/n^{\alpha/\delta}$, we can increase the righthand side of
(\ref{newbound}) as
\begin{eqnarray}
\lefteqn{ 2 M \E U
  \, n_0^{\delta+\epsilon-1}\, |r|
u^{-\delta-\epsilon} \sum_k|\xi_k|  \; [  1_{R_{1,k,n}} +
1_{R_{2,k,n}} \sum_j M |\xi_j|^\kappa  |r|^{\kappa}u^\kappa
1_{\{s\leq t_j\}}  } \label{onceki}
\\  & \;\;\;\;\;\;\;\;\;\;\;\nonumber
\;\;\;\;\;\;\;\;\;\; \displaystyle{
 + 2 \frac{t_k-s}{u/n^{1-\alpha/\delta}}
 \, 1_{R_{3,k,n}}+  \left(2   \frac{t_k}{u/n^{1-\alpha/\delta}}   +
\frac{|s| t_k}{u^2/n^{2(1-\alpha/\delta)}} \right) 1_{R_{4,k,n}} ]}
\end{eqnarray}
which is integrable over $0<u<1$  in view of  Appendix B.

We can now use the dominated convergence theorem. Note that
\begin{eqnarray} \label{extra}
 \lim_{n\rightarrow \infty}\frac{1}{n^{1-\alpha/\delta}}\,\partial_u
g(s,u/n^{1-\alpha/\delta},n^{1-\alpha/\delta}r)\!\! &=& \!\!
\lim_{n\rightarrow \infty}
\partial_u [g(s,u/n^{1-\alpha/\delta},n^{1-\alpha/\delta}r)] \\
\!\! &=&  \!\! \nonumber
\partial_u \lim_{n\rightarrow \infty}g(s,u/n^{1-\alpha/\delta},n^{1-\alpha/\delta}r)
\end{eqnarray}
where $g(s,u/n^{1-\alpha/\delta},n^{1-\alpha/\delta}r)$ is given by
\begin{equation} \label{gsuon}
 e^{i\sum_{k=1}^m\xi_k r u\left[f\left(\frac{t_k-s}{u/n^{1-\alpha/\delta}}\right)-
f\left(\frac{-s}{u/n^{1-\alpha/\delta}}\right)\right]} -1
-i\sum_{k=1}^m\xi_k r u\!\! \left[f\!
\left(\!\frac{t_k-s}{u/n^{1-\alpha/\delta}} \! \right)  -  f \!
\left( \! \frac{-s}{u/n^{1-\alpha/\delta}}\! \right) \right]
\end{equation}
and we have
\begin{equation} \label{limituon} \lim_{n\rightarrow \infty}
\left[f\left(\frac{t_k-s}{u/n^{1-\alpha/\delta}}\right)-
f\left(\frac{-s}{u/n^{1-\alpha/\delta}}\right)\right] =
1_{\{0<s<t_k\}}(s)
\end{equation}
 To see (\ref{limituon}), one takes the limit in regions $R_{1,k},
\ldots , R_{4,k}$, separately. Fig.3 illustrates the function
$\tilde{f}(\cdot):=f(\frac{\cdot\, -\, s}{u/ n^{1-\alpha/\delta}})$
over these regions where we consider $\tilde{f}(t)-\tilde{f}(0)$ as
$n \rightarrow \infty$. By (\ref{kucuklimit}), (\ref{extra}),
(\ref{gsuon}) and (\ref{limituon}), we take the limit of
(\ref{byparts3}) and then revert  the integration by parts to get
the limiting characteristic function as
\[
 \exp \int_{-\infty}^{\infty}\int_0^{\infty}
\int_{-\infty}^{\infty}  [ e^{i\sum_{k=1}^m\xi_k r u
1_{\{0<s<t_k\}}}-1 -i\sum_{k=1}^m\xi_k r u 1_{\{0<s<t_k\}}
 ]\: \lambda u^{-\delta -1}
 ds\,   du\,\gamma(dr) \; .
\]
But, this is the characteristic function of $\sum_{k=1}^m\xi_k
Z(t_k)$, where
\[
Z(t)=\int_{-\infty}^{\infty}\int_0^{\infty} \int_{-\infty}^{\infty}
r\, u \, 1_{\{0<s<t\}} \tilde{N}'(ds,du,dr)
\]
for a Poisson random measure $N'$ with mean measure $\mu'=\lambda ds
u^{-\delta -1}du \gamma(dr)$, and $\tilde{N}'=N'-\mu'$. This
characterizes the limiting process by Lemma \ref{lemma3}.

To complete the proof of weak convergence, it is sufficient to show
that $\mathbb{E}|Z_n(t)-\mathbb{E}Z_n(t)|^{1+\kappa}\leq C t^b$ for
some $b>1$ and $C>0$ in view of the proof of Theorem \ref{ICR}. In
the present theorem, we need a finer estimate given in \cite[Lemma
2]{bahr} and used in \cite[Lemma 6]{kaj}. We have
\begin{equation} \label{chrbd}
\mathbb{E}|Z_n(t)-\mathbb{E}Z_n(t)|^{1+\kappa} \leq a\int_0^{\infty}
\left(1-e^{-2 I_n} \right) \xi^{-2-\kappa} d\xi
\end{equation}
where
\[
I_n= \int \int \int \left(1-\cos\left(\xi \,n^{1-\alpha/\delta}r u
\left[f\left(\frac{t-s}{u}\right) -
f\left(\frac{-s}{u}\right)\right]\right) \right)\mu_n(ds,du,dr)
\]
and $a=(\int_0^{\infty}(1-\cos x)x^{-2-\kappa}dx)^{-1}$, which is
finite with $0<\kappa\leq 1$. Substituting $\mu_n$ and applying
integration by parts, we get
\begin{equation} \label{ay}
I_n=\int \int \int \partial_u k(s,u,n^{1-\alpha/\delta}r) \,
\p\{U>nu\} \frac{n^{\alpha}}{h(n^{\alpha/\delta})}\: \lambda\,
 ds\,  du \,\gamma(dr)
\end{equation}
where
\begin{equation}  \label{k}
k(s,u,r)=  1-\cos\left(\xi \,r u \left[f\left(\frac{t-s}{u}\right) -
f\left(\frac{-s}{u}\right)\right]\right)\; .
\end{equation}
For latter use, the partial derivative of $k$ in $u$ is found as
\begin{eqnarray*}
\lefteqn{\partial_u k(s,u,r) =  \sin \left(\xi ur
\left[f\left(\frac{t-s}{u}\right) -
f\left(\frac{-s}{u}\right)\right] \right) }\\
 & \displaystyle{ \cdot \left[ \xi r \left[ f\left(\frac{t-s}{u}\right) -
f\left(\frac{-s}{u}\right)\right] + \xi r \left[
f'\left(\frac{t-s}{u}\right) \, \frac{t-s}{u} -
f'\left(\frac{-s}{u}\right)\,\frac{-s}{u}\right] \right]}
\end{eqnarray*}
Making a change of variable $u$ to $u/n^{1-\alpha/\delta}$ in
(\ref{ay}), we get
\[
 I_n=\int \int
\int \frac{1}{n^{1-\alpha/\delta}}\, \partial_u
k(s,u/n^{1-\alpha/\delta},n^{1-\alpha/\delta}r) \,
\p\{U>n^{\alpha/\delta}u\}\,
\frac{n^{\alpha}}{h(n^{\alpha/\delta})}\: \lambda \, ds \, du \,
\gamma(dr)\: .
\]
Note the similarity of $I_n$ to (\ref{byparts3}). Moreover, the
inequality $\sin x \leq 2^{1-\kappa}|x|^\kappa \wedge 2$ holds since
$\sin x = [e^{ix}-1 + (e^{ix}+1)]/2$ leading to estimates as in
(\ref{and}) and (\ref{and2}).  It follows that
\[
   4M(C+\epsilon) \xi \!\! \int \!\! \int_{0}^\infty  \!\!\!\!
\int |r| \, B(s,u,t,n) \max (u^{-\epsilon},u^{\epsilon})\,
 u^{-\delta}\, \lambda  \, ds \, du \, \gamma(dr)
\]
is an upper bound to $|I_n|$ when it is  evaluated over $u\geq
n_0/n^{\alpha/\delta}$ where $\epsilon$ and $n_0$ are as above and
\begin{eqnarray*}
\lefteqn{ B(s,u,t,n)\leq   1_{R_{1,k,n}}+\left( 1\wedge   M
|\xi|^\kappa |r|^\kappa u^\kappa  \right).\, 1_{R_{2,k,n}}} \\
&  \displaystyle{+2\, \frac{t-s}{u/n^{1-\alpha/\delta}}\:
1_{R_{3,k,n}}+ \left(2\frac{t}{u/n^{1-\alpha/\delta}} +
\frac{|s|t}{u^2/n^{2(1-\alpha/\delta)}}\right)1_{R_{4,k,n}}}
\end{eqnarray*}
 For evaluating $|I_n|$ for smaller values of $u$, we have a bound  similar
to (\ref{onceki}). Therefore, $I_n$ is bounded by an integrable
function uniformly over $n$ in view of the analogous computations in
Appendix B. By dominated convergence theorem, let $I=\lim_{n}I_n$.
We find that
\[
I=\int_{-\infty}^{\infty}\int_0^{\infty} \int_{-\infty}^{\infty}
[1-\cos(\xi ur1_{\{0<s<t\}}) ]\, u^{-\delta-1}\lambda \, ds \,  du
\, \gamma(dr)
\]
by using the same approach for taking the limit of the
characteristic function of the finite dimensional distributions
above. Then, we can write
\begin{equation}  \label{nerdeyseson}
\mathbb{E}|Z_n(t)-\mathbb{E}Z_n(t)|^{1+\kappa} \leq a
\int_0^{\infty} \left(1-e^{-4 I} \right) \xi^{-2-\kappa} d\xi
\end{equation}
for sufficiently large $n$, by (\ref{chrbd}), since $1-e^{-x}$ is
increasing in $x$. Simplifying $I$ further, we have
\[
I=\lambda t \int_{-\infty}^{\infty}\int_0^{\infty}
  [1-\cos(\xi ur ) ]\,
u^{-\delta-1}  du \, \gamma(dr) = \lambda \,t \,\xi^{\delta}
\mathbb{E} |R|^{\delta} \int_0^{\infty}
  (1-\cos u )\,
u^{-\delta-1}  du
\]
where the second equality follows by a change of variable $u$ to
$u/(\xi |r|)$. Define the constant $\tilde{C}$ so that $I=:\tilde{C}
t\xi^{\delta}$. Now, substituting $I$ in (\ref{nerdeyseson}) and
changing $\xi$ to $\xi/t^{1/\delta}$, we get
\[
\mathbb{E}|Z_n(t)-\mathbb{E}Z_n(t)|^{1+\kappa} \leq a \,
t^{\frac{1+\kappa}{\delta}} \int_0^{\infty} \left(1-e^{-4 \tilde{C}
\xi^{\delta}} \right) \xi^{-2-\kappa} d\xi
\]
which concludes the proof as $(1+\kappa)/\delta>1 $. \finish

\begin{figure}
\centering \epsfig{file=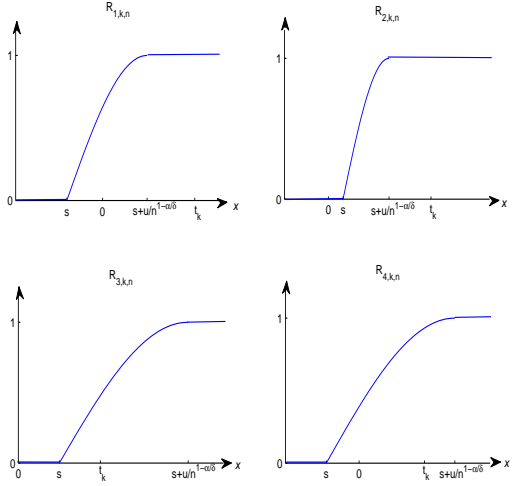, height=2.9in, width=3.2in}
\begin{center} \caption{Graphs of $\tilde{f}(x):=f\left(\frac{x\, -\,
s}{u/ n^{1-\alpha/\delta}}\right)$ in different subregions for
$(s,u/n^{1-\alpha/\delta})$. The function $f$ is arbitrary provided
that  its derivative and itself are Lipschitz continuous almost
everywhere.}
\end{center}
\end{figure}

\begin{rmk}  Note that the stable process obtained in the limit is stable
with a skewness parameter that depends on the distribution of the
rate $R$. Moreover, it has stationary and independent increments.
Therefore, it is also a $\delta$-stable L\'{e}vy motion
\cite[Def.7.5.1]{samor}, but with   scale parameter
\[
\sigma \lambda^{1/\delta}(C_1+C_2)^{1/\delta}
\]
and   skewness parameter
\[
\beta = \frac{C_1-C_2}{C_1+C_2}
\]
by \cite[pg.s 10,11]{samor}, where $C_1 =\mathbb{E}
R^{\delta}1_{\{R>0\}} $ and $C_2 =\mathbb{E} |R|^{\delta}1_{\{R<0\}}
$, see also \cite[pg.217]{bert}. In the context of supply and
demand, one can interpret $C_1$ as the skewness caused by demand and
$C_2$ by supply since they are expected to increase and decrease the
price, respectively.
\end{rmk}

\begin{rmk}
  The weak convergence result given in Theorem \ref{SCR} is proved with
  Skorohod's $J_1$ topology. In \cite{kaj}, the analogous result
  based on (\ref{increasingf}) has been omitted. The convergence is shown
  with $M_1$ topology instead of $J_1$ in \cite{resnick} where
  the  effect function is assumed to be monotone increasing in the
  context of workload input to the system. The authors heuristically
    argue  that some of the individual loads are too large
    \cite[Rmk.4.2]{resnick}.
  On the other hand,   \cite{jedidi} proves weak convergence with
  $M_1$ topology considering that the limit process has jumps.
  However, $J_1$ topology also works as shown above. The interplay
  between $J_1$ and $M_1$ is discussed in \cite{avram} for sums of
  moving averages. It is proved that $J_1$ convergence cannot hold
 because adjacent jumps of this process can coalesce in the limit.
  An intuitive explanation is
given  as the jump of the limiting process occurring  from a
staircase of several jumps.  Under certain conditions, $M_1$
convergence is shown instead.  We have a simpler situation where
each arrival of the scaled process generates a jump of the limit
process as evident from (\ref{limituon}).
\end{rmk}

The following theorem is based on the simpler form of the mean
measure.

\begin{thm} \label{SCRMandelbrot} Suppose the function $f:\mathbb{R}\rightarrow \mathbb{R}$ is Lipschitz
continuous with $f(x)=0$ for all $x\leq 0$, $f(x)=f(1)=1$ for all
$x\geq 1$ and is also differentiable with $f'$ satisfying  a
Lipschitz condition a.e., and $\mathbb{E}|R|^{1+\kappa}<\infty$ for
some $0<\kappa\leq 1$ with $1+\kappa>\delta$. Let
\[
\tilde{Z}_n(t)=\int_{-\infty}^{\infty}\int_0^{\infty}
\int_{-\infty}^{\infty} nr\,u\left[f\left(\frac{t-s}{u}\right)-
f\left(\frac{-s}{u}\right)\right]\tilde{N}_n(ds,du,dr)
\]
where $\tilde{N}=N-\mu$ and
\[ \mu_n(ds,du,dr)=  \lambda \,
n^{-\delta} u^{-\delta-1}\, ds\, du\, \gamma(dr).
\]
 Then, the process
  $\{\tilde{Z}_n(t), t\geq 0\}$, for $1<\delta<2$, converges in law to
  \[
(\lambda \,\mathbb{E} R^{\delta}1_{\{R>0\}})^{1/\delta}\, L_1(t) +
(\lambda \,\mathbb{E} |R|^{\delta}1_{\{R<0\}})^{1/\delta} \,
L_2(t)
\]
 as $n\rightarrow\infty$, where $L_1$ and $L_2$
are  independent $\delta$-stable L\'{e}vy motions with mean 0, and
skewness intensity $1$ and $-1$, respectively.
\end{thm}

\prf We will give only a sketch of the proof due to its similarities
with the previous theorem. The characteristic function for the
finite dimensional distributions of $\tilde{Z}$ can be written as
\[
  \exp \int_{-\infty}^{\infty}\int_0^{\infty}
\int_{-\infty}^{\infty} g(s,u,nr )\, \lambda\,
n^{-\delta}\:u^{-\delta -1}
 ds\,   du\,\gamma(dr)
\]
with $g$  as  in (\ref{g}). Making a change of variable $u$ to
$u/n$, we get
\begin{equation} \label{sonchr}
  \exp \int_{-\infty}^{\infty}\int_0^{\infty}
\int_{-\infty}^{\infty} g(s,u/n,nr )\, \lambda \:u^{-\delta -1}
 ds\,   du\,\gamma(dr)
\end{equation}
Now, $g(s,u/n,nr)$ is similar to (\ref{gsuon}) and we take a
similar limit to (\ref{limituon}) with $n^{1-\gamma/\delta}$
replaced by $n$. This is justified by dominated convergence
theorem
 since the integrand in (\ref{sonchr}) can be bounded as in the
proof of Theorem \ref{SCR}. Convergence in $D(0,\infty)$ follows
along the same lines, this time with $k(s,u/n,nr)$ in $k$ of
(\ref{k}).
 \finish

The simpler form of scalings in Theorems 3 and 5 facilitate  neat
interpretations in terms of the parameters of the price process. In
Theorem 3, $\lambda$ is scaled as $n^2 \lambda$ and $r$ is scaled as
$r/n$, which means that the trading occurs  more frequently, but in
smaller quantities and yields a fractional Brownian motion limit. In
contrast, a stable process is obtained if the rate of trading
decreases while its effect rate increases since $\lambda$ is scaled
as $ \lambda/n^\delta$ and $r$ is scaled as $nr$ in Theorem 5.



\newpage

\noindent {  \bf Appendix A}

 We show that the right hand side of (\ref{bound}) is finite.
  When the right hand
side of (\ref{bound}) is splitted over different regions,  checking
the finiteness of the integrals over  $R_{1,k},R_{3,k},R_{4,k}$
reduces to showing that
\begin{eqnarray*}
 \int_{-\infty}^{0}\int_{-s}^{t_k -s}  u^{-\delta}
\max(u^{-\epsilon},u^\epsilon)\,du \, ds  + \int_0^{t_k}\int_{t_k
-s}^\infty (t_k-s) u^{-\delta-1} \max(u^{-\epsilon},u^\epsilon)\,
du\,  ds \\  + \int_{-\infty}^{0}\int_{t_k -s}^\infty
\left(\frac{1}{u}+\frac{|s|}{u^2}\right) u^{-\delta}
\max(u^{-\epsilon},u^\epsilon)\,du \, ds
\end{eqnarray*}
is finite. This is indeed true when we choose $\epsilon>0$ such that
\begin{equation} \label{ep2}
1<\delta -\epsilon<\delta <\delta+\epsilon<2 \: .
\end{equation}
In region $R_{2,k}$, we have
\begin{equation} \label{I1}
I:=\int  \int_0^{t_k}\int_{0}^{t_k-s} |r|[ 1\wedge \sum_j M
|\xi_j|^\kappa |r|^\kappa
 u^{\kappa}1_{\{s\leq t_j\}} ]u^{-\delta}
\max(u^{-\epsilon},u^\epsilon)\,du \, ds\, \gamma(dr)
\end{equation}
If $t_j>t_k$, it can be observed from Fig.2 that the integral
reduces to that over region $R_{2,k}$. If $t_j<t_k$, then the
integral over $R_{2,k}$ yields an upper bound. That is, we can
replace $1_{\{s\leq t_j\}}$ by the constant function $1$ and get
\begin{eqnarray}
\lefteqn{I \leq \label{I}
       \E |R|^{1+\kappa}\sum_j |\xi_j|^\kappa \,
     \int_{0}^{\bar{u}}  \int_0^{t_k}
   u^{\kappa -\delta}\max
(u^{-\epsilon},u^{\epsilon}) \,   ds \, du }\\ &  \qquad
 + \, \displaystyle{\E |R| \int_{\bar{u}}^{t_k}  \int_0^{t_k}
  u^{-\delta}\, \max
(u^{-\epsilon},u^{\epsilon})\,   ds \, du} \nonumber
\end{eqnarray}
where $\bar{u}$ denotes a cutoff value of $u$ such that $\sum_j
|\xi_j|^\kappa u^\kappa$ is too large in (\ref{I}), and we use the
fact that $t_k-u\leq t_k$ for $u\geq 0$ after changing the order of
integration for $u$ and $s$ in (\ref{I1}). Then, the right hand side
of (\ref{I}) is finite if we choose $\epsilon>0$ such that
\[
1<\delta -\epsilon<\delta <\delta+\epsilon<1+\kappa
\]
which clearly satisfies  (\ref{ep2}) since $\kappa\leq 1$.

\newpage

\noindent {  \bf Appendix B}

In this part, we show that (\ref{42}) is integrable with respect to
$\max(u^{-\epsilon},u^\epsilon)u^{-\delta}\,du\, ds $. Substituting
the limits of integration in regions $R_{1,k,n},R_{3,k,n},R_{4,k,n}$
shown by $I_1,I_3,I_4$, respectively, we have
\begin{eqnarray*}
I_1 &=&
\int_{-\infty}^0\int_{-\tilde{n}s}^{\tilde{n}t_k-\tilde{n}s}\max(u^{-\epsilon},u^\epsilon)
\, u^{-\delta}\, du \, ds  \\
I_3&=& \int_{0}^{t_k} \int_{\tilde{n}t_k-\tilde{n}s}^{\infty}
(\tilde{n}t_k-\tilde{n}s) \max(u^{-\epsilon},u^\epsilon)
\, u^{-\delta-1}\, du \, ds\\
I_4&=& \int_{-\infty}^0 \int_{\tilde{n} t_k-\tilde{n}s}^{\infty}
\left(\frac{2}{u/\tilde{n}}+\frac{|s|}{u^2/\tilde{n}^2}\right)\max(u^{-\epsilon},u^\epsilon)
\, u^{-\delta}\, du \, ds
\end{eqnarray*}
where we put $\tilde{n}=n^{1-\alpha/\delta}$. The integrals
$I_1,I_3,I_4$ are finite for
\[
1<\delta -\epsilon<\delta <\delta+\epsilon<2
\]
since
\begin{eqnarray*}
\int_{-\infty}^0\int_{-\tilde{n}s}^{\tilde{n}t_k-\tilde{n}s}
\, u^{-\tilde{\delta}}\, du \, ds  &=&C_1\frac{t_k^{2-\tilde{\delta}}}{\tilde{n}^{\tilde{\delta}-1}} \leq C_1\, t_k^{2-\tilde{\delta}} \\
  \int_{0}^{t_k} \int_{\tilde{n}t_k-\tilde{n}s}^{\infty} (\tilde{n}t_k-\tilde{n}s)
 u^{-\tilde{\delta}-1}\, du \, ds &=& C_2\frac{t_k^{2-\tilde{\delta}}}{\tilde{n}^{\tilde{\delta}-1}} \leq C_2\, t_k^{2-\tilde{\delta}} \\
  \int_{-\infty}^0 \int_{\tilde{n} t_k-\tilde{n}s}^{\infty}
\left(\frac{2}{u/\tilde{n}}+\frac{|s|}{u^2/\tilde{n}^2}\right)
u^{-\tilde{\delta}}\, du \, ds &\leq & C_3
\frac{t_k^{2-\tilde{\delta}} +t_k^{-\tilde{\delta}}
+t_k^{1-\tilde{\delta}}       }
    {\tilde{n}^{\tilde{\delta}-1}} \\& \leq & C_3 (t_k^{2-\tilde{\delta}}+
t_k^{-\tilde{\delta}} +t_k^{1-\tilde{\delta}})
\end{eqnarray*}
for $1<\tilde{\delta}<2$ and $\tilde{n}\geq 1$, where
$C_1,C_2,C_3\in \mathbb{R}$. In $R_{2,k,n}$, we have
\[
I_2=\int_0^{t_k}\int_0^{\tilde{n}t-\tilde{n}s} \left(1\wedge \sum_j
M |\xi_j|^\kappa |r|^\kappa u^\kappa 1_{\{s\leq t_j \}}\right)
\max(u^{-\epsilon},u^\epsilon)\, u^{-\delta}\, du\, ds
\]
As in Appendix A, we consider two intervals $[0,\bar{u}]$ and
$(\bar{u},t_k]$ to evaluate this integral. Over the first interval,
it is finite for $1<\tilde{\delta}<1+\kappa$, and over the latter,
it is proportional to $\tilde{n}^{1-\tilde{\delta}}$ which is
bounded by 1. As a result, (\ref{42}) is finite if we choose
$\epsilon>0$ such that
\[
1<\delta -\epsilon<\delta <\delta+\epsilon<1+\kappa \;.
\]


\begin{thebibliography}{9}

\bibitem{akcay} Z. Ak\c{c}ay, M. \c{C}a\u{g}lar, An Agent Based Stock Price Model, SPA
2007, Illinois, USA.

\bibitem{avram} F. Avram, M. Taqqu, Weak Convergence of Sums of
Moving Averages in the $\alpha$-Stable Domain of Attraction, The
Annals of Probability, 20 (1992), 483-503.

\bibitem{bayrak} E. Bayraktar , U. Horst, R. Sircar,
A limit Theorem for Financial Markets with Inert Investors,
Mathematics of Operations Research, 31 (2006), 789-810.

\bibitem{bert} J. Bertoin, L'{e}vy Processes, 2007,
 Cambridge University Press.

\bibitem{bill} P. Billingsley, Convergence of Probability
Measures, 1999, Second Edition,  John Wiley.

\bibitem{bingham} N.H. Bingham, C.M. Goldie, J.L. Teugels, Regular
Variation, 1987, Cambridge University Press.

\bibitem {demand4} G.-I. Bischi, M. Gallegati, L. Gardini, R.
Leombruni, A. Palestrini, Herd behavior and nonfundamental asset
price fluctuations in financial markets, Macroeconomic Dynamics, 10
(2006), 502-528.

\bibitem{demand5} D. Challet, A. Chessa, M. Marsili and Y.-C.
Zhang, From Minority Games to real markets, Quantitative Finance
Volume 1 (2001), 168-176.

\bibitem{mandel} R. Cioczek-Georges, B. B. Mandelbrot,
Alternative Micropulses and Fractional Brownian Motion, Stochastic
Processes and their Applications, 64 (1996), 143-152.

\bibitem{caglar} M. \c{C}a\u{g}lar,
A Long-Range Dependent Workload Model for Packet Data Traffic.
Mathematics of Operations Research, 29 (2004), 92-105.

\bibitem{world} M. \c{C}a\u{g}lar, Weak Limits of Infinite Source Poisson
Pulses, 7th World Congress in Probability and Statistics, 2008,
Singapore.

\bibitem{cinlarlec1} E. \c{C}{\i}nlar, Lecture Notes in Probability
Theory, Princeton University, 2003. Probability and Stochastics,
Springer (to appear, 2010).

\bibitem{demand1} G. Iori, A microsimilation of traders activity
in the stock market: the role of heterogeneity, agents' interactions
and trade frictions, Journal of Economic Behavior and Organization,
49 (2002), 269-285.

\bibitem{chartist-demand1} J. D. Farmer, S. Joshi, The Price dynamics of
common trading strategies, Journal of Economic Behavior and
Organization, 49 (2002), 149-171.

\bibitem{fasen} V. Fasen, G. Samorodnitsky A Fluid Cluster Poisson
Process Can Look Like a Fractional Brownian Motion Even in the Slow
Growth Regime, Advances in Applied Probability, 41 (2009), 393-427.

\bibitem{feller} W. Feller, An Introduction to Probability Theory
and its Applications, Vol.2, 2nd edn, 1971, Wiley, New
York.

\bibitem{complex} W. Fulks, Complex Variables, Marcel Dekker,
Inc., 1993.

\bibitem{gaigalas} R. Gaigalas, I. Kaj, Convergence to Scaled
Renewal Processes and a Packet Arrival Model, 9 (2003), 671-703.

\bibitem{gaig} R. Gaigalas, A Poisson Bridge between Fractional Brownian Motion
and Stable L\'{e}vy Motion, Stochastic Processes and their
Applications, 116 (2006), 447-462.

\bibitem {demand3} F. Ghoulmie, R. Cont and J.-P. Nadal,
Heterogeneity and feedback in an agent-based market model, Journal
of Physics: Condensed Matter, 17 (2005), S1259-S1268.

\bibitem {chartist-demand2} I. Giardina, J-P. Bouchaud, M. M\'{e}zard,
Microscopic models for long ranged volatility correlations, Physica
A, 299 (2001), 28-39.

\bibitem{jedidi} W. Jedidi, J. Almhana, V. Choulakian, R.  McGorman,
 The Poisson Shot Noise Traffic Model, and
Functional Convergence to Stable Processes, (2010), Preprint.

\bibitem {chartist-demand3} X-Z. He, F. H. Westerhoff, Commodity markets,
price limiters and speculative price dynamics, Journal of Economic
Dynamics and Control, 29 (2005), 1577-1596.

\bibitem{kall} O. Kallenberg Foundations of Modern Probability,
1997, Springer.

\bibitem{kaj} I. Kaj and M. S. Taqqu, Convergence to fractional Brownian motion and
to the Telecom process: the integral representation approach,
Brazilian Probability School, 10th anniversary volume, Eds. M.E.
Vares, V. Sidoravicius, 2007, Birkhauser.

\bibitem{kingman} J.F.C. Kingman, Poisson Processes, 1993, Clarendon Press,
Oxford.

\bibitem{scharf} C. Kl\"{u}ppelberg, T. Mikosch and A. Scharf,
Regular Variation in the Mean and Stable Limits for Poisson Shot
Noise, Bernoulli, 9 (2003), 467-496.

\bibitem{klup} C. Kl\"{u}ppelberg and C. K\"{u}hn, Fractional
Brownian motion as a weak limit of Poisson shot noise processes-
with applications to  finance, Stochastic Processes and their
Applications, 113 (2004), 333-351.

\bibitem{konstan} T. Konstantopoulos and S. Lin . Macroscopic Models for
Long-Range Dependent Network Traffic, Queueing Systems, 28 (1998)
215-243.

\bibitem{kurtz} T. Kurtz . Limit Theorems for Workload Input Models. F.P.
Kelly, S. Zachary, I. Ziedins, Eds.,  Stochastic Networks: Theory
and Applications,  1996, Clarendon Press, Oxford.

\bibitem{mikres} T. Mikosch, S. Resnick, H. Rootz\'{e}n and A. Stegeman .
Is Network Traffic Approximated by Stable L\'{e}vy Motion or
Fractional Brownian Motion?   Annals  of Applied Probability,  12
(2002) 23-68.

\bibitem{mikosch} T. Mikosch and G. Samorodnitsky,
Scaling Limits for Cumulative Input Processes, Mathematics of
Operations Research, 32 (2007), 890-918.

\bibitem[Samorodnitsky and Taqqu(1994)]{samor}
G. Samorodnitsky and M.S. Taqqu  Stable Non-Gaussian Random
Processes, 1994, Chapman \& Hall.

\bibitem{demand2} K. S. Sznajd-Weron, R. Weron, A simple model
of price formation, International Journal of Modern Physics C, 13
(2002), 115-123.

\bibitem{pipiras} V. Pipiras, M.S. Taqqu, J.B. Levy, Slow, Fast and
Arbitrary Growth Conditions for Renewal-Reward Processes When Both
the Renewals and the Rewards Are Heavy-Tailed, Bernoulli, 10 (2004),
121-163.

\bibitem{resnick} S. Resnick, E. van den Berg, Weak Convergence of
High-Speed Network Traffic Models, J. Appl. Prob., 37 (2000),
575-597.

\bibitem{sato} K. Sato, L\'{e}vy Processes and Infinitely Divisible
Distributions, 1999, Cambridge University Press.

\bibitem{bahr} B. von Bahr, C.G. Esseen, Inequalities for the
\emph{r}th Absolute Moment of a Sum of Random Variables, $1\leq
r\leq 2$, Ann. Math. Stat., 36 (1965), 299-303.

\end{thebibliography}
\end{document}